\documentclass{article}
\usepackage{amsmath}
\usepackage{amssymb}
\usepackage{float}
\usepackage[top=2.54cm, bottom=2.54cm, left=2.5cm, right=2.5cm]{geometry}
\usepackage{graphicx}
\usepackage{natbib}
\usepackage{caption}
\usepackage{subcaption}

\usepackage{fancyhdr}
\begin{document}

\title{The Arab Spring:\\ A Simple Compartmental Model for the Dynamics of a Revolution}
\author{John Lang and Hans De Sterck\\\small{Department of Applied Mathematics, University of Waterloo}}
\date{\today}
\maketitle

\begin{abstract}
	The self-immolation of Mohamed Bouazizi on December 17, 2011 in the small Tunisian city of Sidi Bouzid, set off a sequence of events culminating in the revolutions of the Arab Spring.  It is widely believed that the Internet and social media played a critical role in the growth and success of protests that led to the downfall of the regimes in Egypt and Tunisia.  However, the precise mechanisms by which these new media affected the course of events remains unclear.  We introduce a simple compartmental model for the dynamics of a revolution in a dictatorial regime such as Tunisia or Egypt which takes into account the role of the Internet and social media.  An elementary mathematical analysis of the model identifies four main parameter regions: stable police state, meta-stable police state, unstable police state, and failed state.  We illustrate how these regions capture, at least qualitatively, a wide range of scenarios observed in the context of revolutionary movements by considering the revolutions in Tunisia and Egypt, as well as the situation in Iran, China, and Somalia, as case studies.  We pose four questions about the dynamics of the Arab Spring revolutions and formulate answers informed by the model.  We conclude with some possible directions for future work.
\end{abstract}

\section{Introduction and Motivation}
\label{sec:Motive}

\begin{quotation}
	``After decades of political stagnation... new winds of hope were felt in the Middle East, accompanied by a new catchword making the rounds in the American media, `Arab Spring'...  The age of the old patriarchs, it appeared, was nearing its end.  And the new media - satellite television, mobile phones, the Internet - were often regarded as having precipitated this development by undermining governments' hegemonic control over the flow of information.''
\end{quotation}

\paragraph{}  When Albrecht Hofheinz wrote these words he was referring to modest advancements being made in democracy and political liberalization in a handful of Middle Eastern countries in 2005 \citep{Hofheinz05}.  He did not foresee the events sparked by Mohamed Bouazizi's self-immolation on December 17, 2010 that ultimately led to the Arab Spring revolutions.  Nevertheless, his analysis of new media and their impact on Arab society are eerily prescient, especially considering that in 2005 social media was either in its infancy or completely non-existent\footnote{Facebook was launched in 2004 and was still an invitation-only service in 2005 \citep{Phillips07}.  Youtube was founded in early 2005 \citep{Youtube} and Twitter was not founded until the spring of 2006 \citep{Picard11}.}.  He concludes,

\begin{quotation}
	``The Internet is one factor that in tandem with others (satellite TV, youth culture, and the `globalization' of consumer products, social networks, and ideational configurations) is creating a dynamic of change that is helping to erode the legitimacy of traditional authority structures in family, society, culture/religion, and also the state, and thus creating pressure for reform.''
\end{quotation}

\paragraph{}  Consistent with \citet{Hofheinz05}, the predominant view is that the Internet and social media played a critical role in the Arab Spring of 2010-2011 \citep{HowardEtAl11, LotanEtAl11, Alterman11, BBC11, KhamisVaughn11, Pollock11, Saletan11, Shirky11, Stepanova11, ZhuoEtAl11}.  Although most of these opinions are based on anecdotal evidence, some rigorous work has been done attempting to determine a link between social media and protests using, for example, Twitter data \citep{HowardEtAl11, LotanEtAl11}.  Needless to say, many questions remain, for example:

\begin{enumerate}
	\item  How can a small number of active social media users and relatively low Internet penetration\footnote{According to \citet{HowardEtAl11} approximately 25\% of Tunisians and 10\% of Egyptians had used the Internet at least once prior to the Arab Spring.} have a dramatic effect on the stability of a regime?,
	\item  How is it that regimes manage to seem so stable until the revolution is underway?,
	\item  Why did the January 28 - February 1, 2011, Internet shutdown in Egypt not have a greater inhibitory effect on protests?, and
	\item  Why is it that some regimes fall in a matter of weeks, others fight to a stalemate, and still others survive relatively unscathed?
\end{enumerate}

\paragraph{}  Answering these four questions, among others, using a consistent and unified approach is not trivial given the complexity of the situations in Arab Spring countries.  Adding to the difficulty of such a task is the impossibility of running counter-factual experiments to verify conclusions.  However, the goal of this paper is to show how mathematical modelling can be useful in this situation.  It should be noted that models of opinion/norm formation \citep{CentolaEtAl05}, conflict \citep{AtkinsonEtAl11,Kress12}, and even revolution \citep{Kuran91} already exist.  However, these models either do not apply specifically to the Arab Spring revolutions or are highly complex.  Furthermore, although complex models may in principle be able to offer a more complete description, they also have limitations.  More detailed models typically require additional assumptions and the calibration of a large number of parameters.  This makes complex models analytically intractable, difficult to interpret, and computationally expensive to simulate.  In this paper we attempt to create a simple model that is nevertheless able to capture essential features of Arab-Spring-type revolutions, and also can be used to explore possible answers to the questions raised above.

\paragraph{}  We develop our model in Section \ref{sec:Intro} and provide an elementary mathematical analysis in Section \ref{sec:Analysis}.  This is followed by Section \ref{sec:Interp} which expands on the interpretation of our model by considering various case studies.  This section also explains how our model can help to find answers to the four questions posed above.  Section \ref{sec:Tunis} first describes the Arab Spring events in Tunisia, applying our model with parameters that are fixed for the lifetime of the revolution.  Next we consider the Egyptian revolution in Section \ref{sec:Egypt}, where we allow model parameters to evolve over the course of the revolution in order to incorporate external influences.  Finally, Section \ref{sec:Other} briefly discusses how our model can be applied to several other states, including Iran, China, and Somalia.  The paper is concluded in Section \ref{sec:Conc} with a summary of our findings and a discussion of future work.

\section{Simple Compartmental Model}
\label{sec:Intro}
\paragraph{}  In order to facilitate the development, interpretation, and analysis of our model we begin by stating it and defining the terminology used in Section \ref{sec:StateModel}.  Section \ref{sec:InterpJustModel} then provides the interpretation and justification of each model term and parameter introduced in Section \ref{sec:StateModel}.  

\subsection{Statement of the Model}
\label{sec:StateModel}
\paragraph{}  The function $r(t)$ represents the fraction of \emph{protesters} or \emph{revolutionaries} in the population at time $t$.  The model which we use to describe the dynamics of the revolutions in Tunisia and Egypt is given by a single differential equation for $r(t)$,

\begin{equation}
	\label{eq:ODE}
	\dot{r} = \underbrace{c_1\mbox{ } v(r;\alpha)\mbox{ } (1-r)}_{g(r)} - \underbrace{c_2\mbox{ } p(r;\beta)\mbox{ } r}_{d(r)},
\end{equation}

where parameters $\alpha,\beta\in(0,1)$ and $c_1,c_2>0$, where $\dot{x}$ denotes the time derivative of $x$, and where the functions $g,d:[0,1]\rightarrow\mathbb{R}^+$ are called the \emph{growth} and \emph{decay} terms, respectively, since they model the growth and decay of the fraction of protesters.

\begin{figure}[h]
	\centering
	\begin{subfigure}[t]{3in}
		\centering
		\includegraphics[width=3in]{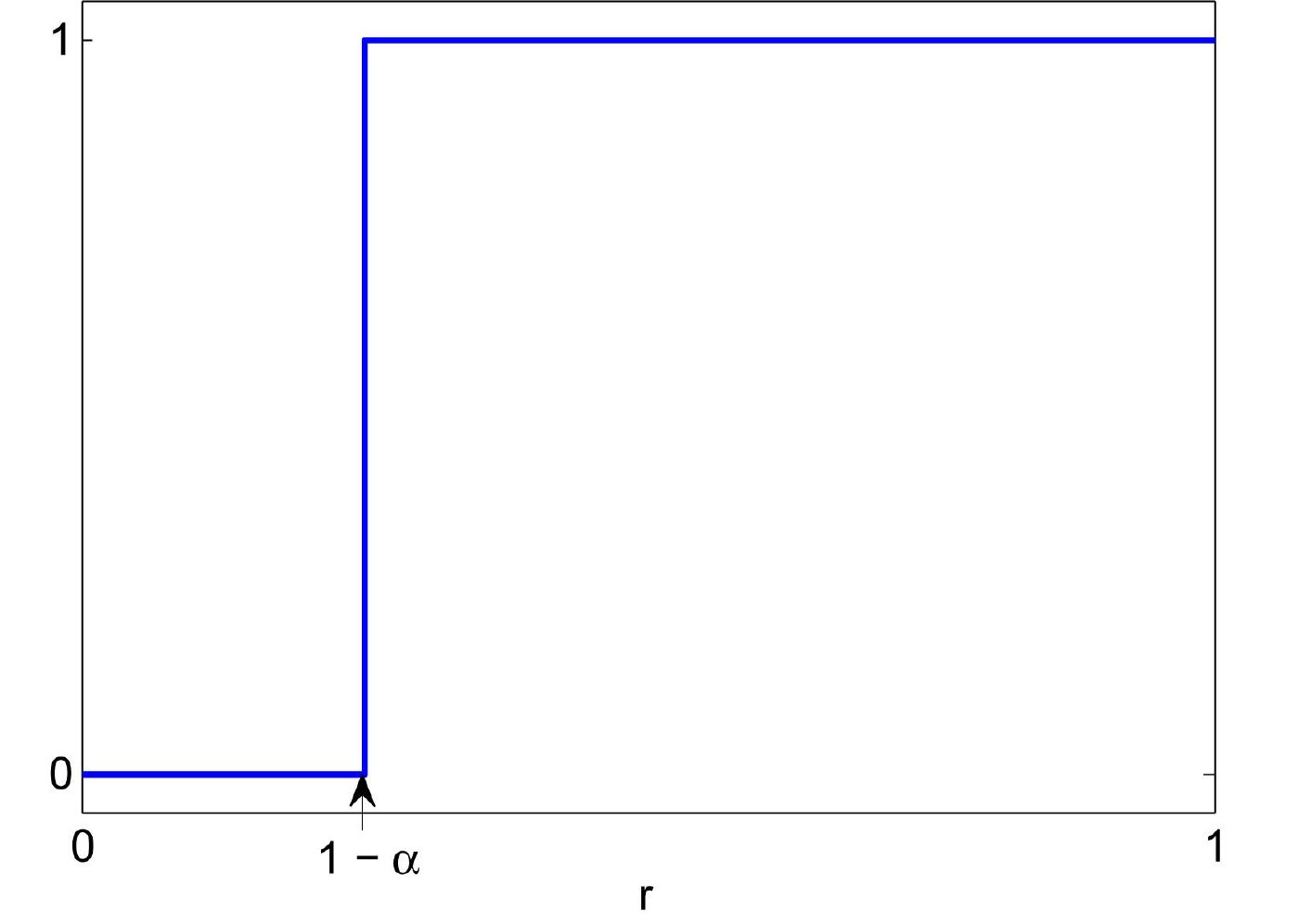}
		\caption{Visibility term, $v(r;\alpha)$}
	\end{subfigure}
	\hfill
	\begin{subfigure}[t]{3in}
		\centering
		\includegraphics[width=3in]{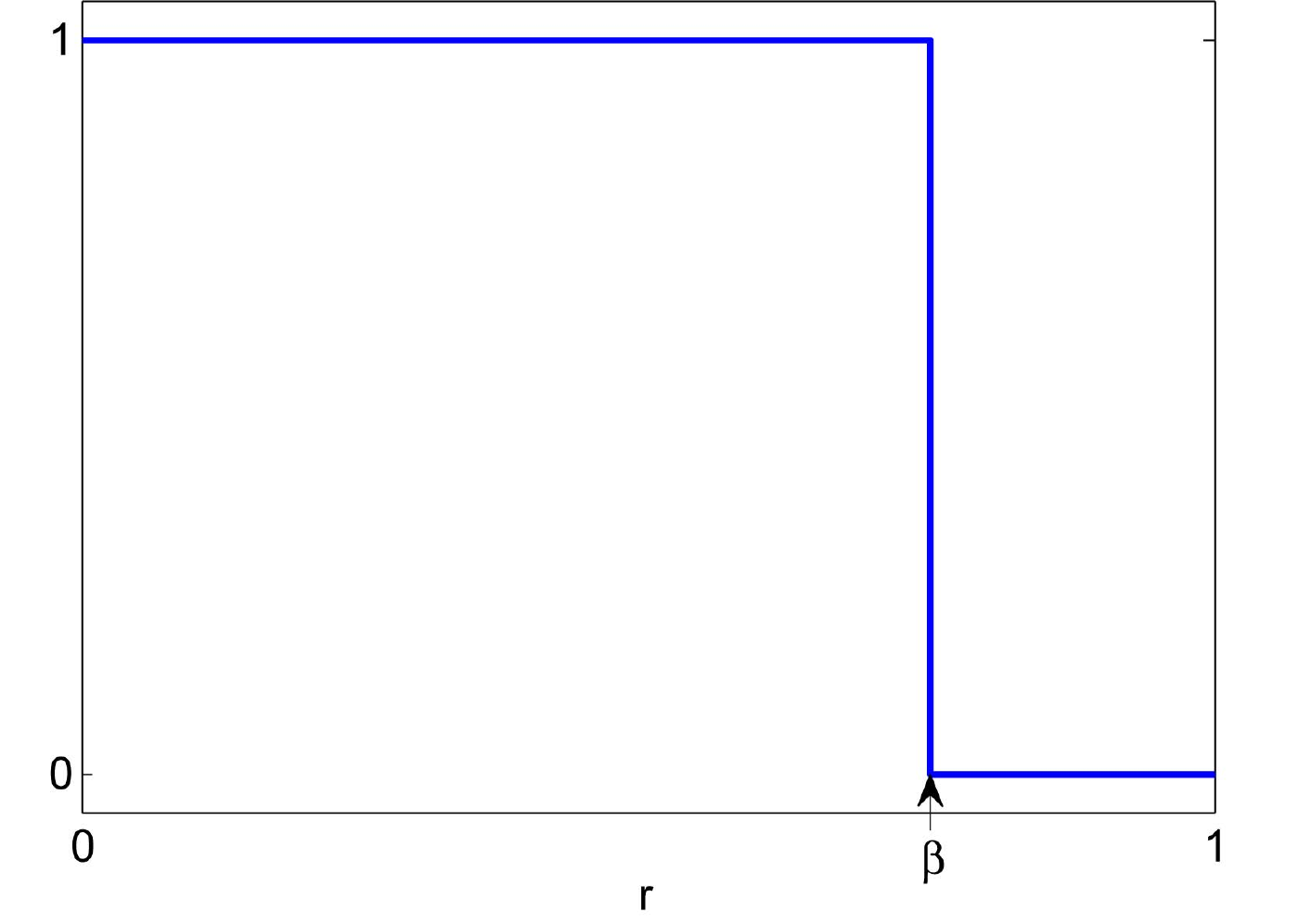}
		\caption{Policing term, $p(r;\beta)$}
	\end{subfigure}
	\caption{Visibility ($v(r;\alpha)$) and policing ($p(r;\beta)$) terms.}
	\label{fig:VP}
\end{figure}

\paragraph{}  Subject to the \emph{visibility term}

$$
	v(r;\alpha) = \left\{\begin{array}{ll} 1 & \mbox{if } r > 1 - \alpha\\ 0 & \mbox{otherwise}.
		\end{array}\right.
$$

the growth term is proportional to $(1-r)$.  We call the proportionality constant, $c_1$, and the parameter, $\alpha$, the \emph{enthusiasm} and \emph{visibility} of the protesters, respectively.  The visibility term is modelled as a step function, which shuts off the growth term when the fraction of protesters is below the \emph{visibility threshold} $1-\alpha$: in our model the fraction of protesters can only grow when the protest movement is sufficiently large to be visible to the general population.  Similarly, subject to the \emph{policing term}

$$
	p(r;\beta) = \left\{\begin{array}{ll} 1 & \mbox{if } r < \beta\\ 0 & \mbox{otherwise}.
		\end{array}\right.
$$

the decay term is proportional to $r$.  We call the proportionality constant, $c_2$, and the threshold parameter, $\beta$, the \emph{policing capacity} and \emph{policing efficiency}, respectively.  The policing term is also modelled as a step function: it shuts down the decay term when the fraction of protesters is above the policing capacity threshold, $\beta$.  The visibility and policing terms are illustrated in Figure \ref{fig:VP}.

\paragraph{}  We observe from equation \eqref{eq:ODE} and Figure \ref{fig:VP} that if $r=0$ or $r=1$ then $\dot{r}=0$, regardless of the values chosen for the parameters.  We say that $r=0$ and $r=1$ are the equilibria of \emph{total state control} and of the \emph{realized revolution}, respectively.  In what follows it will be shown how this model can describe the dynamics of a revolutionary transition from a small initial group of protesters ($r\approx0$) to a full-blown revolution ($r\approx1$).  This is a simple model for revolutionary transitions, parametrized by the four parameters $\alpha$, $\beta$, $c_1$, and $c_2$.

\paragraph{}  Now that we have established the terminology for our model we are able to provide a detailed interpretation and justification in Section \ref{sec:InterpJustModel}.  This will allow us to proceed with Sections \ref{sec:Analysis} and \ref{sec:Interp}, which provide the mathematical analysis of our model and application of the model to various case studies, respectively.

\subsection{Interpretation and Justification of the Model}
\label{sec:InterpJustModel}

\begin{figure}[h]
	\centering
	\includegraphics[width=4in]{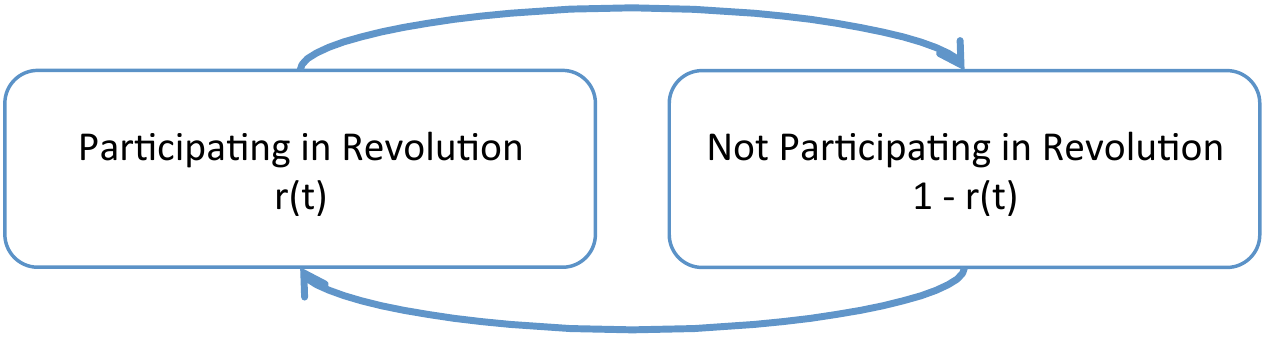}
	\caption{Simple compartmental model for the dynamics of a revolution.}
	\label{fig:Compartments}
\end{figure}

\paragraph{}   Here we explain and justify the model of \eqref{eq:ODE} which describes the process by which citizens engage in revolution, with the specific goal of gaining insight into the effect of enhanced communications technologies \citep{HowardEtAl11, LotanEtAl11}.  In order to arrive at a simple model some simplifying assumptions are necessary.  First, our model is developed for describing rapid revolutionary transitions on a short time scale (of the order of months), and neglects demographic and other long-term effects.  Second, we assume that the regime is very unpopular and that all individuals would privately like to see the regime changed.  The second assumption allows us to divide the population into two compartments: the population participating/not participating in the revolution.  From the first assumption, the sum of both compartments is a constant and the dynamics of one compartment completely determines the dynamics of the other, so it suffices to consider a one-compartment model, see Figure \ref{fig:Compartments}.  We choose to keep track of the population participating in the revolution as a fraction of the total population, $r(t)$.  Note that the fraction of the population available to join the revolution at time $t$ is $1-r(t)$, by the second assumption.  We also note that the first assumption is applied again in Sections \ref{sec:Tunis}-\ref{sec:Egypt} to identify reasonable values for the protester's enthusiasm ($c_1$) and the regime's policing efficiency ($c_2$).  It remains to justify our choice of functional form for $g(r)$ and $d(r)$ in \eqref{eq:ODE}.

\paragraph{}  We assume that the regime is capable of arresting/dispersing protesters at a rate proportional to the size of the revolution, $r$, provided that the number of protesters does not exceed the regime's finite policing capacity, $\beta$.  Provided that no new protesters join the revolution ($v=0$) and that the number of protesters does not exceed the regime's policing capacity ($p=1$), this corresponds to exponential decay in the number of protesters with the timescale determined by the policing efficiency, $c_2$.  We make the further simplifying assumption that the regime loses all ability to punish protesters once the number of protesters exceeds the regime's policing capacity.  These assumptions determine the form of $p(r;\beta)$ as a switching function.

\paragraph{}  Dictatorial regimes are known to keep tight control on the flow of political information through state control of the media and through censorship, and for a good reason \citep{CentolaEtAl05,Dunn11,BBC11,KhamisVaughn11,Kuran91}: if political protests are kept hidden from the general population, protest movements have little chance of growing.  We model this effect by the visibility term, $v(r;\alpha)$, and we make the simple assumption that the visibility term can be modelled as a step function.  As soon as the fraction of protesters reaches the visibility threshold, $1-\alpha$, and is large enough to be visible to the general population, the revolution is assumed to grow with a growth term proportional to $1-r$.  Note that we call $\alpha$ the visibility and $1-\alpha$ the visibility threshold: for large visibility, e.g. $\alpha=0.96$, the visibility threshold is low, $1-\alpha=0.04$, so the general population will become aware of the political protest movement as soon as it has spread to 4\% of the population.  One of the goals of this paper is to investigate how the increased presence of Internet, social media, satellite television, and cell phones may enhance the spread of revolutionary movements.  Indeed, these effects may significantly loosen the control of the regime over the flow of politically sensitive information, and in our model the influence of new media can be taken into account by an increased visibility parameter, $\alpha$.  Note also that, provided the revolution is visible ($v=1$) and exceeds the policing capacity ($p=0$), the growth rate is proportional to $1-r$ with the timescale determined by the protesters' enthusiasm, $c_1$.

\paragraph{}  As a secondary motivation for the step-function form of $v(r;\alpha)$ one can also consider the following.  Given the policing limitations of the regime, the decision of individuals whether or not to act is a collective action problem \citep{Kuran91}.  Thus, the case can be made that the most important factor for individuals deciding to join a revolution is the \emph{perceived size} of the revolution.  If individuals perceive participation in a revolution to be below a certain threshold they will refuse to join the revolution and risk punishment, despite their desire to see the regime fall.  Conversely, above this threshold an individual's desire to see the regime fall overpowers their fear of government reprisal.

\paragraph{}  Due to the simplicity of our four-parameter model it is unable to capture singular one-time events such as the self-immolation of Mohamed Bouazizi on December 17, 2010.  Although these types of events could be modelled stochastically, to keep our model as simple as possible we introduce the concept of \emph{shocks}.  A shock is an event external to our model which nevertheless has an effect on the fraction of revolutionaries ($r$) either directly, or indirectly via a change in the parameters $\alpha$, $\beta$, $c_1$, or $c_2$.  Specifically, we consider direct shocks that produce an instantaneous jump in $r$, denoted $\Delta r$, and indirect shocks that trigger an instantaneous or continuous change in one or more parameters of the model.  A shock of this type to $\alpha$, for example, would be specified by defining $\alpha(t)$.  for simplicity we restrict our attention to instant or linear changes in parameters.

\begin{figure}[h!]
	\centering
	\includegraphics[width=5in]{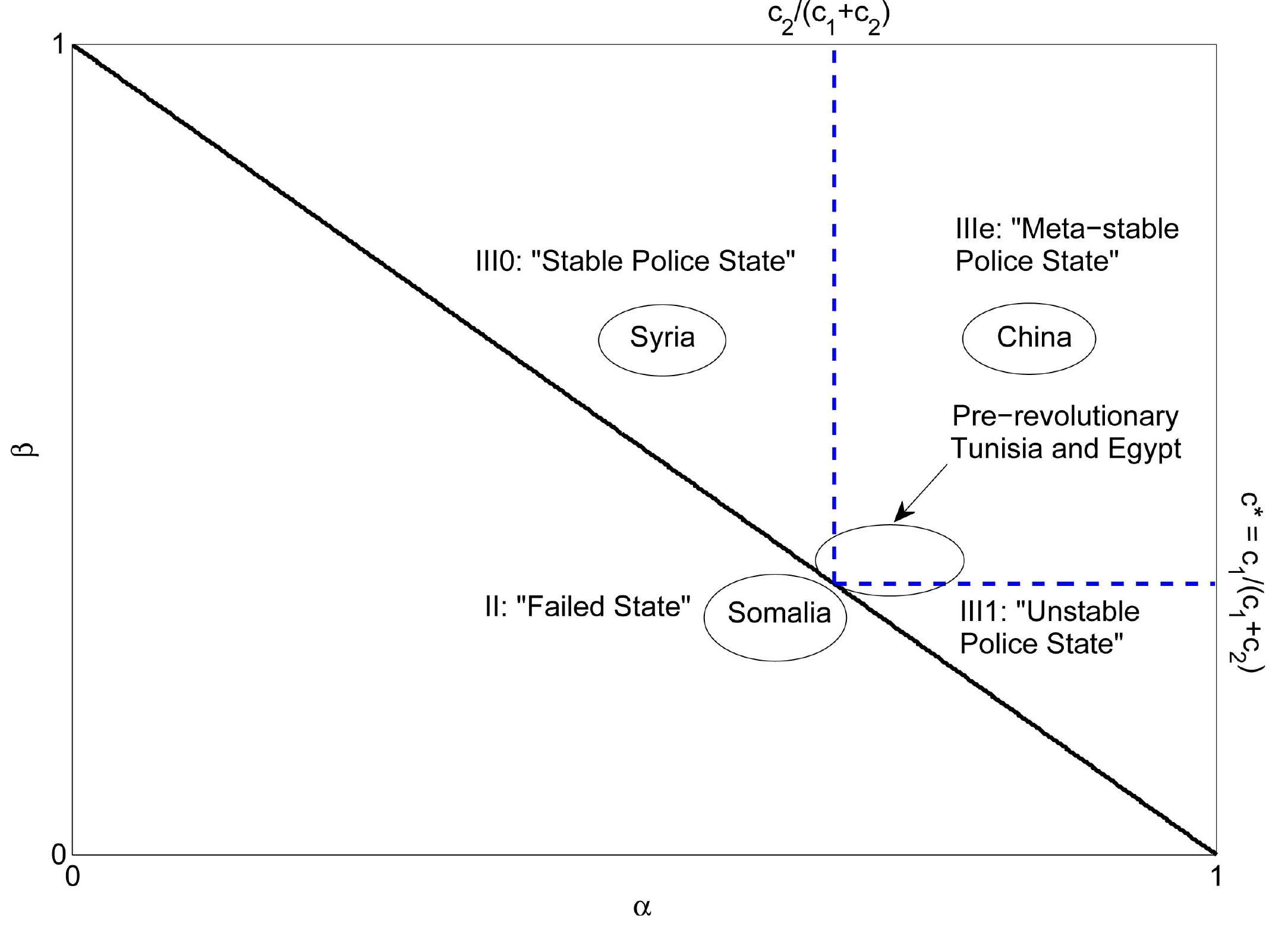}
	\caption{Division of $\alpha-\beta$ parameter space into regions II, IIIe, III0, and III1.  Conceptual summary of case studies of Tunisia, Egypt, Iran, China, and Somalia.}
	\label{fig:Summary}
\end{figure}

\paragraph{}  In Section \ref{sec:Interp} we will discuss how our model can be applied in the case studies of the Tunisian and Egyptian revolutions, and the situations in Iran, China, and Somalia.  Before embarking on this study, however, Section \ref{sec:Analysis} provides a mathematical analysis of our model.  We find that the parameter space of the model can essentially be divided into four regions, which we name regions II, III0, IIIe, and III1 (see Figure \ref{fig:Summary}), and which we will interpret in terms of the dynamic stability of the model solutions in those parameter regions.  In anticipation of the case studies of Section \ref{sec:Interp}, Figure \ref{sec:Analysis} gives a conceptual indication of how the cases of Tunisia, Egypt, Iran, China, and Somalia can be described by our model using parameter choices in specific parts of the parameter space.

\paragraph{}  The assumptions made in developing our model are crude, but as illustrated by the analysis and case studies presented below, they apply sufficiently well to some of the Arab Spring revolutions that they can be used to formulate a simple model that captures some essential features of these revolutions.  In particular, we emphasize that  our model focuses specifically on the types of rapid transitions that have characterized the revolutions in Tunisia and Egypt, and that all individuals would privately like to see the regime change.  It is clear that there are countries for which these assumptions and the resulting model do not apply, as discussed in Section \ref{sec:Conc} on future work.

\section{Elementary Mathematical Analysis of the Model}
\label{sec:Analysis}

\begin{figure}[h]
	\centering
	\includegraphics[width=4in]{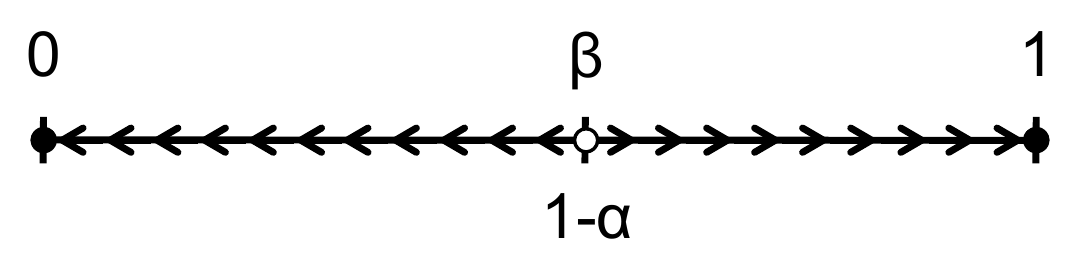}
	\caption{Region I with $\alpha+\beta=1$.  Closed (open) circles represent locally asymptotically stable (unstable) equilibria.  Left (right) arrows indicate regions where $\dot{r}<0$ ($\dot{r}>0$).}
	\label{fig:I}
\end{figure}

\begin{figure}[h!]
	\centering
	\includegraphics[width=4in]{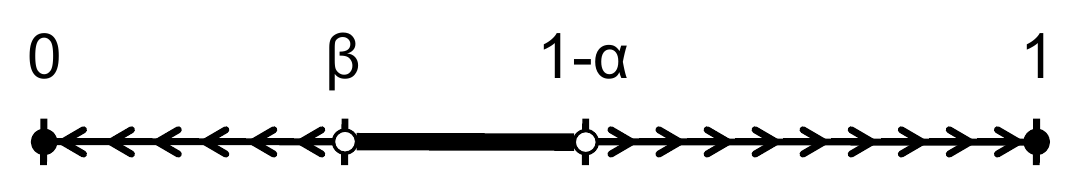}
	\caption{Region II with $\alpha + \beta <1$.  Closed (open) circles represent locally asymptotically stable (unstable) equilibria.  The thick line in $(\beta,1-\alpha)$ indicates stable equilibria.  Left (right) arrows indicate regions where $\dot{r}<0$ ($\dot{r}>0$).}
	\label{fig:II}
\end{figure}

\begin{figure}[h!]
	\centering
	\begin{subfigure}[t]{4in}
		\centering
		\includegraphics[width=4in]{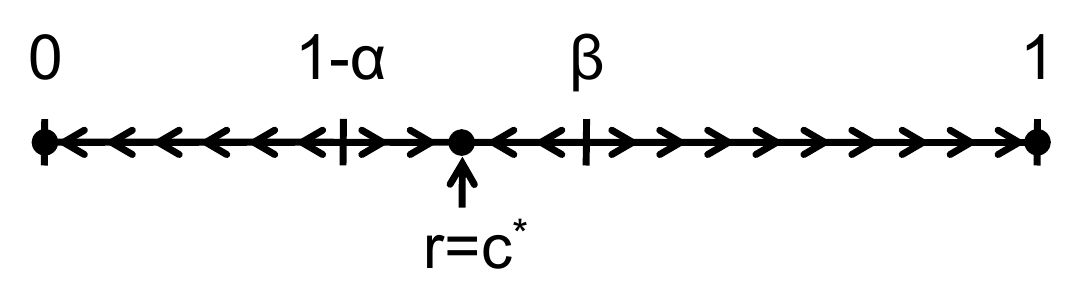}
		\caption{Region IIIe with $\alpha+\beta>1$ and $c^*\in(1-\alpha,\beta)$}
	\end{subfigure}
	\\
	\begin{subfigure}[t]{4in}
		\centering
		\includegraphics[width=4in]{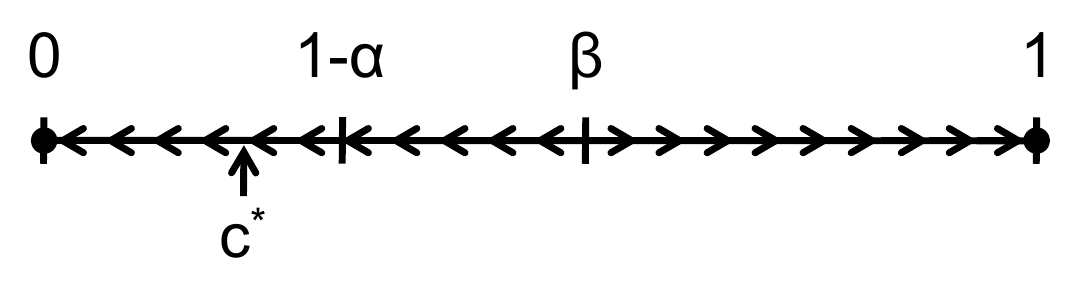}
		\caption{Region III0 with $\alpha+\beta>1$ and $c^*<1-\alpha$}
	\end{subfigure}
	\\
	\begin{subfigure}[t]{4in}
		\centering
		\includegraphics[width=4in]{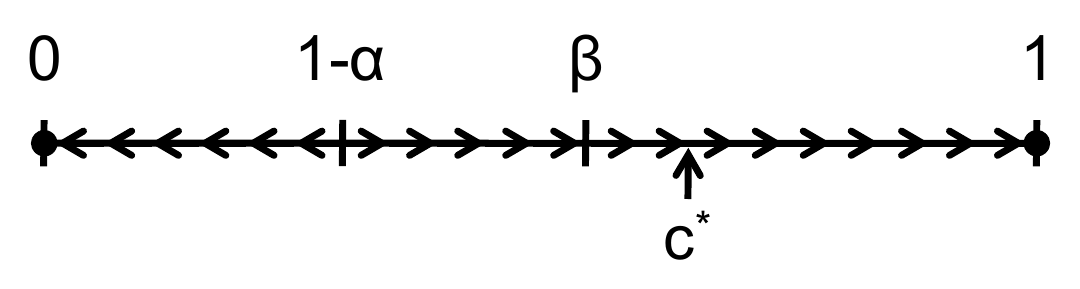}
		\caption{Region III1 with $\alpha+\beta>1$ and $c^*>\beta$}
	\end{subfigure}
	\caption{Regions IIIe, III0 and III1.  Closed (open) circles represent locally asymptotically stable (unstable) equilibria.  Left (right) arrows indicate regions where $\dot{r}<0$ ($\dot{r}>0$).}
	\label{fig:III}
\end{figure}

\paragraph{}  The mathematical analysis of the dynamics of the model given in \eqref{eq:ODE} breaks down into three distinct cases characterized by regions in the $\alpha-\beta$ plane: $\alpha + \beta = 1$, $\alpha + \beta <1$, and $\alpha + \beta > 1$, which are summarized in Figures \ref{fig:I}, \ref{fig:II}, and \ref{fig:III}, respectively.  An interpretation for these cases is given at the end of this section and in subsequent sections.

\paragraph{Region I: $\alpha+\beta=1$}  
\paragraph{}  When $r<\beta=1-\alpha$ we have $v(r;\alpha)=0$ and $p(r;\beta)=1$, so $r=0$ is a locally asymptotically stable equilibrium with basin of attraction $(0,\beta)$.  Similarly, when $r>\beta=1-\alpha$ we have $v(r;\alpha)=1$ and $p(r;\beta)=0$, so $r=1$ is also a locally asymptotically stable equilibrium with basin of attraction $(1-\alpha,1)$.  Finally, because $v(1-\alpha;\alpha)=p(\beta;\beta)=0$ it follows that $r=\beta=1-\alpha$ is a locally unstable equilibrium.

\paragraph{Region II: $\alpha + \beta<1$}
\paragraph{}  As above, $r=0$ and $r=1$ are locally asymptotically stable equilibria with basins of attraction $(0,\beta)$ and $(1-\alpha,1)$, respectively.  When $r\in[\beta,1-\alpha]$ we have $v(r;\alpha)=p(r;\beta)=0$, so all $r\in(\beta,1-\alpha)$ are locally stable equilibria and $r\in\{\beta,1-\alpha\}$ are unstable equilibria.

\paragraph{Region III: $\alpha + \beta >1$}
\paragraph{}  Analogously to the previous two cases, we have the locally asymptotically stable equilibria $r=0$ and $r=1$ with basins of attraction $(0,1-\alpha]$ and $[\beta,1)$, respectively.  Restricting our attention to the interval $(1-\alpha,\beta)$ and solving the algebraic equation $\dot{r}=0$ gives $r = \frac{c_1}{c_1+c_2}$.  We define $c^*=\frac{c_1}{c_1+c_2}$ and observe that our analysis breaks down into a further three sub-cases.

\paragraph{Region IIIe:}  If $\frac{c_1}{c_1+c_2}\in(1-\alpha,\beta)$ then there exists a third equilibrium $r=\frac{c_1}{c_1+c_2}$ and this equilibrium is locally asymptotically stable with basin of attraction $(1-\alpha,\beta)$.

\paragraph{Region III0:}  If $\frac{c_1}{c_1+c_2}<1-\alpha$ then the region $(1-\alpha,\beta)$ lies in the basin of attraction of $r=0$.

\paragraph{Region III1:}  If $\frac{c_1}{c_1+c_2}>\beta$ then the region $(1-\alpha,\beta)$ lies in the basin of attraction of $r=1$.

\paragraph{}  The above results are summarized in Figure \ref{fig:Summary}, which illustrates the relationship between regions II, IIIe, III0, and III1 in the $\alpha-\beta$ plane.  Region I, because it is one-dimensional, is unlikely to manifest itself and so we mostly disregard it in what follows.  For the remaining regions we now introduce terminology to ease future discussion.  States with parameters in region II have uncountably many stable equilibria.  These equilibria occur because the policing capacity of the regime is too low to clear the protesters and the visibility is too low to attract new protesters.  We therefore interpret region II as corresponding to a \emph{failed state}.  Regions III0, IIIe, and III1 differ only in the stability of the interval $(1-\alpha,\beta)$.  For region III0 the interval $(1-\alpha,\beta)$ lies in the basin of attraction of total state control ($r=0$).  Because of the contribution of $(1-\alpha,\beta)$ to the stability of the regime, we refer to region III0 as a \emph{stable police state}.  Analogously, we refer to region III1 as an \emph{unstable police state}.  Region IIIe introduces an intermediate state, $r=c^*$, which lies between the equilibria of total state control ($r=0$) and of the realized revolution ($r=1$).  We therefore refer to $r=c^*$ as the equilibrium of \emph{civil unrest} and to region IIIe as a \emph{meta-stable police state}.  

\paragraph{}  Solutions to the model with parameters in the failed state (region II) are relatively straightforward, because for any value of $r\in[0,1]$ at most one of the visibility ($v(r;\alpha)$) or policing ($p(r;\beta)$) terms is ``switched on''.  Behaviour of solutions with parameters in the stable police state, meta-stable police state, and unstable police state (regions III0, IIIe, and III1, respectively) are less obvious, because for parameters in these regions we have $v(r;\alpha)=1$ and $p(r;\beta)=1$ whenever $r\in(1-\alpha, \beta)$.  Therefore, in order to facilitate the interpretation of our model, the case studies of Section \ref{sec:Interp} will provide sample time traces of solutions to our model for parameters in regions III0, IIIe, and III1.

\section{Interpretation}
\label{sec:Interp}
\paragraph{}  Section \ref{sec:Analysis} provided a mathematical analysis for our model.  However, before we proceed with the case studies of Tunisia and Egypt it will be useful to identify some reasonable values for the protesters' enthusiasm ($c_1$) and the policing efficiency of the regime ($c_2$).  Since we have assumed that the revolutions occur as rapid transitions over a on a short time scale (on the order of months), we take as our guide the observed time scales in the revolutions we want to model.  In the absence of government repression ($p=0$) and with visibility ($v=1$) we assume that the revolution would spread to 90\% of the population within one month for the types of revolutions we want to study\footnote{There were 29 days between Mohamed Bouazizi's self-immolation on December 17, 2010, and Ben Ali's resignation in Tunisia on January 14, 2011 and 18 days between the January 25, 2011, Tahrir Square protests and Mubarak's resignation in Egypt on February 11, 2011 \citep{HowardEtAl11, BBC11, Alterman11}.  Of course, there is no way to determine (a) the exact start date for the revolutions in Egypt or Tunisia, (b) how many people had joined these revolutions by the fall of the regimes, or (c) how these revolutions would have proceeded in the absence of government repression.  Nevertheless, spread of support of the revolution to 90\% of the population in one month seems to be at least within the correct order of magnitude.}.  Measuring time in months and solving the equation $\dot{r} = c_1(1-r)$ with conditions $r(0)=0$ and $r(1) = 0.9$ implies 

$$c_1 = \log(10)\approx2.30.$$  

Similarly, in the absence of new revolutionary recruits ($v=0$) and with perfect policing capacity of the regime ($p=1$) we assume that to clear $90\%$ of revolutionaries would take one day\footnote{Again, this is only a crude order-of-magnitude estimate based on the fact that Egyptian forces managed to clear Tahrir Square in approximately 24 hours after the January 25th protests began \citep{BBC11}.}.  Again measuring time in months and solving $\dot{r}=-c_2r$ with $r(0)=r_0$ and $r(\frac{1}{30}) = 0.1\mbox{ }r_0$ implies 

$$c_2 = 30\log(10)\approx  69.1.$$

\paragraph{}  We are now prepared to consider how our model might be applied to the two Arab Spring revolutions of Tunisia (Section \ref{sec:Tunis}) and Egypt (Section \ref{sec:Egypt}).  We also consider how the situations in Iran, China, and Somalia might fit into the framework developed above (Section \ref{sec:Other}).  Sections \ref{sec:Tunis}-\ref{sec:Other} together provide examples of regimes fitting each of the parameter regions of our model: the failed, stable police, meta-stable police, and unstable polices states, as was summarized in Figure \ref{fig:Summary}.  Furthermore, by the end of Section \ref{sec:Other} we will have seen how our model might be applied to address the four questions posed in Section \ref{sec:Motive}.

\subsection{Case Study: Tunisia}
\label{sec:Tunis}
\paragraph{}  The Arab Spring had its first manifestation in the Jasmine Revolution of Tunisia where the Internet, and in particular social media (i.e Twitter, Facebook, Youtube, etc...), is credited as a catalyst facilitating regime change \citep{HowardEtAl11, Pollock11, Stepanova11}.  Below, we apply the above model in an attempt to better understand how social media may have influenced the revolution.

\paragraph{}  Internet, social media, satellite TV, and cell phone communications technologies may empower protesters by enhancing their

\begin{enumerate}
	\renewcommand{\labelenumi}{(\arabic{enumi})}
	\item capacity for organization and coordination \citep{Beckett11, BBC11, Pollock11},
	\item ability to assess the current public support for the revolution \citep{Alterman11, BBC11, KhamisVaughn11, Pollock11, Saletan11, ZhuoEtAl11}, and
	\item awareness of the nature and severity of government repression \citep{BBC11, Schneider11}.
\end{enumerate}

These and related effects enter into our model via the visibility ($\alpha$) and enthusiasm ($c_1$) parameters in the growth term of the model.  Thus, before considering how our model might be applied to the specific circumstances of the Tunisian revolution, we elaborate on the three points brought up above.

\paragraph{}  As we discussed in Section \ref{sec:Intro}, the decision of whether or not to protest is largely a coordination problem \citep{Kuran91}.  If individuals protest individually then the state is capable of severe retaliation, however, if individuals protest in sufficient numbers then the state loses its ability to punish.  This realization lead activists to use the Internet to coordinate the initial protests in Tunis \citep{BBC11, Pollock11}.  Once protests were underway, technologies such as SMS and Twitter messaging were used between co-revolutionaries, for example by communicating which streets were the most/least obstructed by security forces \citep{BBC11}.  This enhances the speed with which revolutionaries mobilize, and in the context of our model this corresponds to an increased enthusiasm of the revolutionaries ($c_1$).  Social media and the Internet also contributed to the relatively leaderless way in which the Tunisian revolution developed.  Contrasting to revolutions with a more hierarchical leadership structure, a leaderless revolution is difficult if not impossible to disrupt by targeting only a handful of individuals \citep{Beckett11, BBC11}.  This also corresponds to an increased $c_1$.  

\paragraph{}  In our model growth of the revolution is subject to the visibility switching term, $v(r;\alpha)$.  The revolution can only grow when it is visible, i.e. when $r>1-\alpha$.  Through censorship dictatorial regimes attempt to control protests by ensuring that they remain virtually invisible to the general population.  In other words, in order to prevent growth of small protests into larger movements the regime will attempt to keep the visibility ($\alpha$) low, and the visibility threshold ($1-\alpha$) high.  As a consequence, the general population remains unaware of, and hence incapable of joining, small protests.  The Internet, social media, satellite TV, and cell phones all work towards increasing $\alpha$ by disrupting the regime's monopoly on the distribution of information.  In Tunisia the Internet and social media created a virtual space where Tunisians could express their true opinions with minimal censorial oversight or fear of reprisal \citep{BBC11, KhamisVaughn11, Pollock11}.  Critically, this new interconnectivity allowed Tunisians to better gage the true level of support for regime change.  Together with cell phones, social media sites vastly sped up the speed with which information travelled, allowing Tunisians - and the entire world - to follow the revolution with unprecedented detail and speed \citep{LotanEtAl11, Saletan11, ZhuoEtAl11}.  Traditonal media lent its credibility to this new wellspring of information by corroborating and then rebroadcasting stories relating to the size of the revolution and the regime's brutal response \citep{Alterman11}.

\paragraph{}  Awareness of the brutality and severity of the government's reaction may increase both the visibility of protesters, $\alpha$, as well as their enthusiasm, $c_1$.  The enthusiasm for the revolution may be directly affected by increasing resentment of the regime.  In contrast, the effect on $\alpha$ is likely to be through a secondary chanel.  Specifically, otherwise apolitical individuals are induced to join the revolution \citep{BBC11, Schneider11}, presumably by lowering their personal thresholds for participation.

\paragraph{}  We now have an idea of how the Internet, social media, satellite TV, and cell phones might influence the parameters of our model, specifically the visibility ($\alpha$) and enthusiasm ($c_1$) parameters in the growth term.  Next we present a simplified timeline of major events during the Tunisian revolution \citep{BlightEtAl12, HowardEtAl11, BBC11,Rifai11}.

\begin{itemize}
	\item December 17, 2010: Mohamed Bouazizi self-immolates in the city of Sidi Bouzid.
	\item December 18, 2010: Protests erupt in Sidi Bouzid.  Protesters begin recording and uploading videos of the protests and police response to the Internet.
	\item December 22, 2010: Houcine Falhi commits suicide by electrocution in the midst of a demonstration in Sidi Bouzid.
	\item December 27 - 28, 2010: Protests break out in the capital, Tunis.  President Ben Ali denounces protests in televised address.
	\item January 5, 2011: Mohamed Bouazizi dies from burn injuries.
	\item January 14 - 15, 2011: Ben Ali resigns and flees to Saudi Arabia.  Interim government formed.
	\item January 17 - 24, 2011: Protests continue increasing in size to 40,000 - 100,000 individuals.
	\item Februrary 27, 2011: Protest of at least 100,000 Tunisians forces Prime Minister Mohamed Ghannouchi to resign.
\end{itemize}

\begin{figure}[h]
	\centering
	\begin{subfigure}[t]{2.75in}
		\centering
		\includegraphics[width=2.75in]{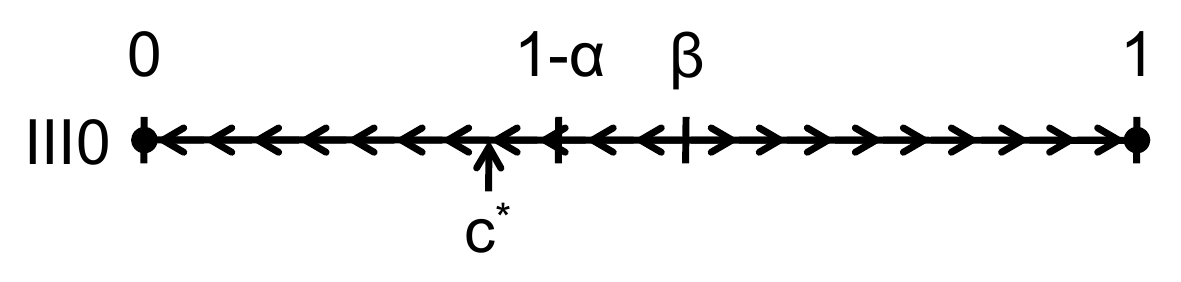}\\
		\includegraphics[width=2.75in]{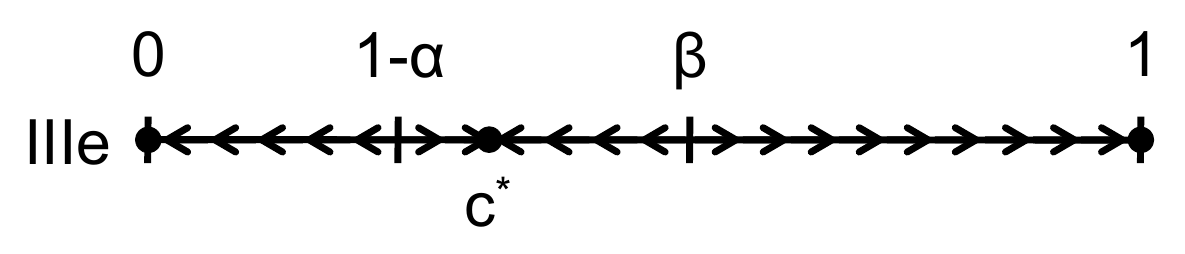}\\
		\includegraphics[width=2.75in]{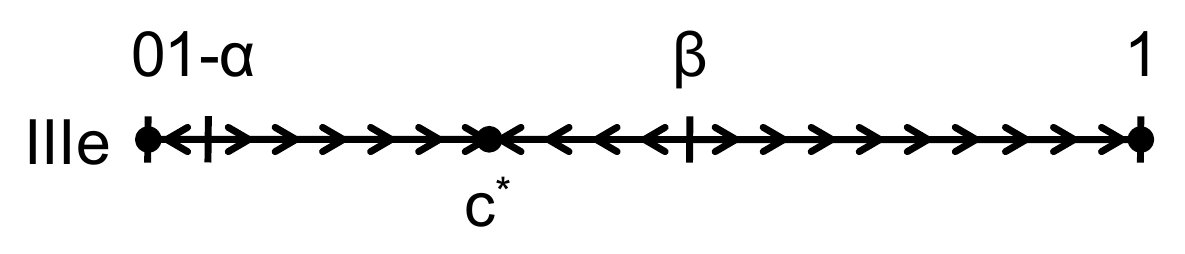}\\
		\caption{Increasing $\alpha$:  While $\alpha<1-c^*$, increasing $\alpha$ has no effect on the basin of attraction of either $r=0$ or $r=1$.  Once $\alpha>1-c^*$, $r=c^*$ becomes a locally asymptotically stable equilibrium, thus reducing the size of a perturbation required to leave the basin of attraction of $r=0$.  Finally, as $\alpha$ becomes large enough the perturbation needed to return to $r=0$ from $r=c^*$ becomes larger than the one required to proceed from $r=c^*$ to $r=1$, and the perturbation required to leave the basin of attraction of $r=0$ becomes increasingly small.}
		\label{fig:TunisiaAlpha}
	\end{subfigure}
	\hspace{5mm}
	\begin{subfigure}[t]{2.75in}
		\centering
		\includegraphics[width=2.75in]{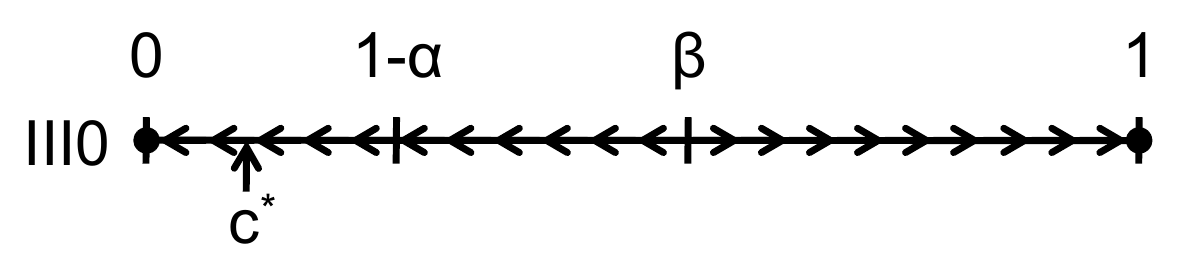}\\
		\includegraphics[width=2.75in]{IIIealpha1}\\
		\includegraphics[width=2.75in]{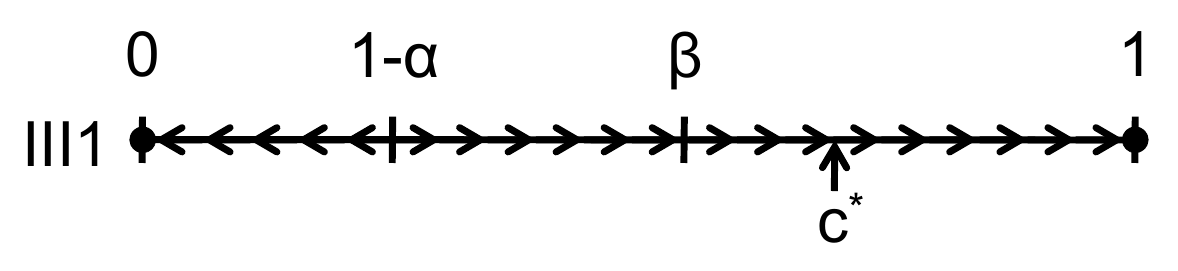}\\
		\caption{Increasing $c_1$:  Increasing $c_1$ has no effect on the basin of attraction of either $r=0$ or $r=1$ until $\alpha>1-c^*$, when $r=c^*$ becomes a locally asymptotically stable equilibrium.  As $c^*$ approaches $\beta$ the equilibrium $r=c^*$ becomes farther from (closer to) to the basin of attraction of $r=0$ ($r=1$).  As a result, the size of a perturbation required to leave $r=c^*$ for $r=0$ ($r=1$) increases (decreases) as $c^*$ increases.  Finally, as $c^*>\beta$ the equilibrium $r=c^*$ disappears and the basin of attraction of $r=1$ becomes $(1-\alpha,1)$. }
		\label{fig:TunisiaC1}
	\end{subfigure}
	\caption{The effect of increasing $\alpha$ and $c_1$ on the existence and stability of equilibria.}
	\label{fig:TunisiaLine}
\end{figure}

\paragraph{}  Our model can now be applied to show how an increase in $\alpha$ and/or $c_1$ might increase the likelihood of a successful revolution.  Given the many decades of stability within Tunisia before the revolution, we assume that the situation pre-Internet is best described by parameters in the stable police state (region III0, $c^*<1-\alpha$).  To see clearly how increasing visibility ($\alpha$) or enthusiasm ($c_1$) of protesters then affects the dynamics of our model we now consider varying these parameters separately at first, starting with $\alpha$.  At first an increasing $\alpha$ leaves the basin of attraction of the total state control equilibrium ($r=0$) unchanged, and hence, has no impact on the overall stability of the regime.  However, once $\alpha > 1 - c^*$ the regime passes from the stable police state (region III0) to the meta-stable police state (region IIIe), the basin of attraction for total state control equilibrium shrinks from $(0,\beta)$ to $(0,1-\alpha]$, and the locally asymptotically stable civil unrest equilibrium ($r=c^*$) is created with basin of attraction $(1-\alpha,\beta)$, see Figure \ref{fig:TunisiaAlpha}.  The size of the perturbation required to leave the basin of attraction of the total state control equilibrium has thus been decreased.  As $\alpha$ continues to increase, the basin of attraction of the civil unrest equilibrium grows.  Depending on the value of $\beta$ and $c^*$, $\alpha$ may be able to grow to the point that a revolution succeeding, i.e. leaving the basin of attraction of the civil unrest equilibrium for that of the realized revolution ($r=1$), requires a smaller perturbation than the revolution being crushed, i.e. leaving the civil unrest equilibrium for the basin of attraction of total state control, see Figure \ref{fig:TunisiaAlpha}.  So, increasing $\alpha$ (a) first moves us into the meta-stable police state (region IIIe), thus decreasing the size of the perturbation needed to leave total state control ($r=0$), and then (b) decreases the relative likelihood of a revolution dying out (returning to $r=0$ from $r=c^*$) as opposed to succeeding (leaving $r=c^*$ for $r=1$).

\paragraph{}  Again starting from a stable police state (region III0), we consider the effect of increasing enthusiasm of protesters ($c_1$), see Figure \ref{fig:TunisiaC1}.  For $c_1<\beta$ the effect of increasing $c_1$ is qualitatively similar to increasing $\alpha$.  Specifically, increasing $c_1$ has no effect on the basin of attraction of total state control ($r=0$) until $c^*>1-\alpha$ ($\alpha>1-c^*$) at which point the regime moves from the stable police state (region III0) to the meta-stable police state (region IIIe), with effects as described above.  As $c_1$ increases with $c^*$ in the region $(1-\alpha,\beta)$, the equilibrium of civil unrest ($r=c^*$) moves farther from the basin of attraction of total state control and closer to the basin of attraction of the realized revolution ($r=1$).  The effects are, again, as above.  Further increasing $c_1$ such that $c^*$ surpasses $\beta$ produces the qualitatively distinct effect of moving from the meta-stable police state (region IIIe) to the unstable police state (region III1).  The equilibrium of civil unrest disappears and the basin of attraction of the realized revolution is extended from $(\beta,1)$ to $(1-\alpha,1)$.

\begin{figure}[h]
	\centering
	\begin{subfigure}{3in}
		\centering
		\includegraphics[width=3in]{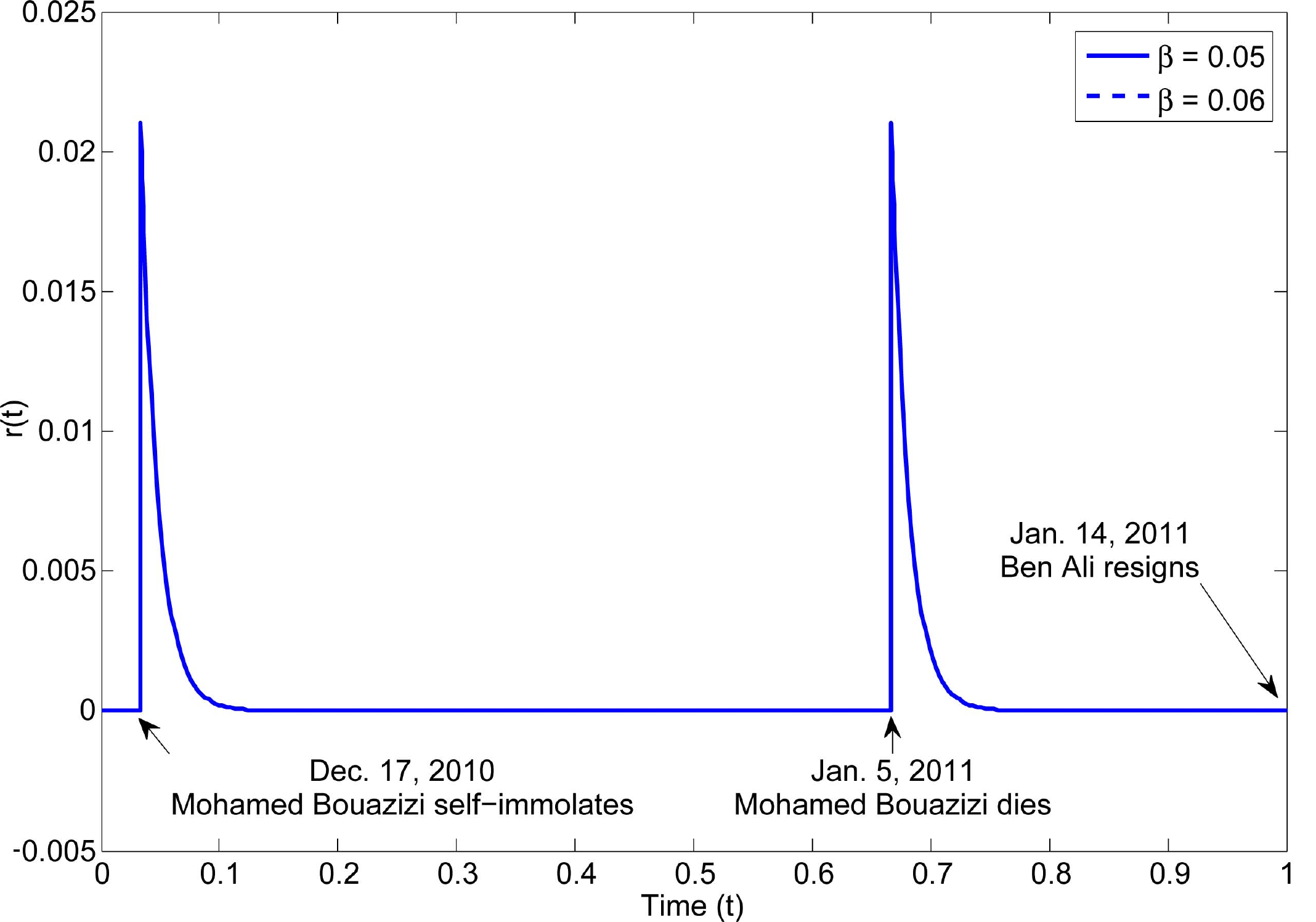}
		\caption{Stable police state (region III0) with $\alpha=0.96$.}
		\label{fig:TunisiaTimeAlphaA}
	\end{subfigure}
	\hfill
	\begin{subfigure}{3in}
		\centering
		\includegraphics[width=3in]{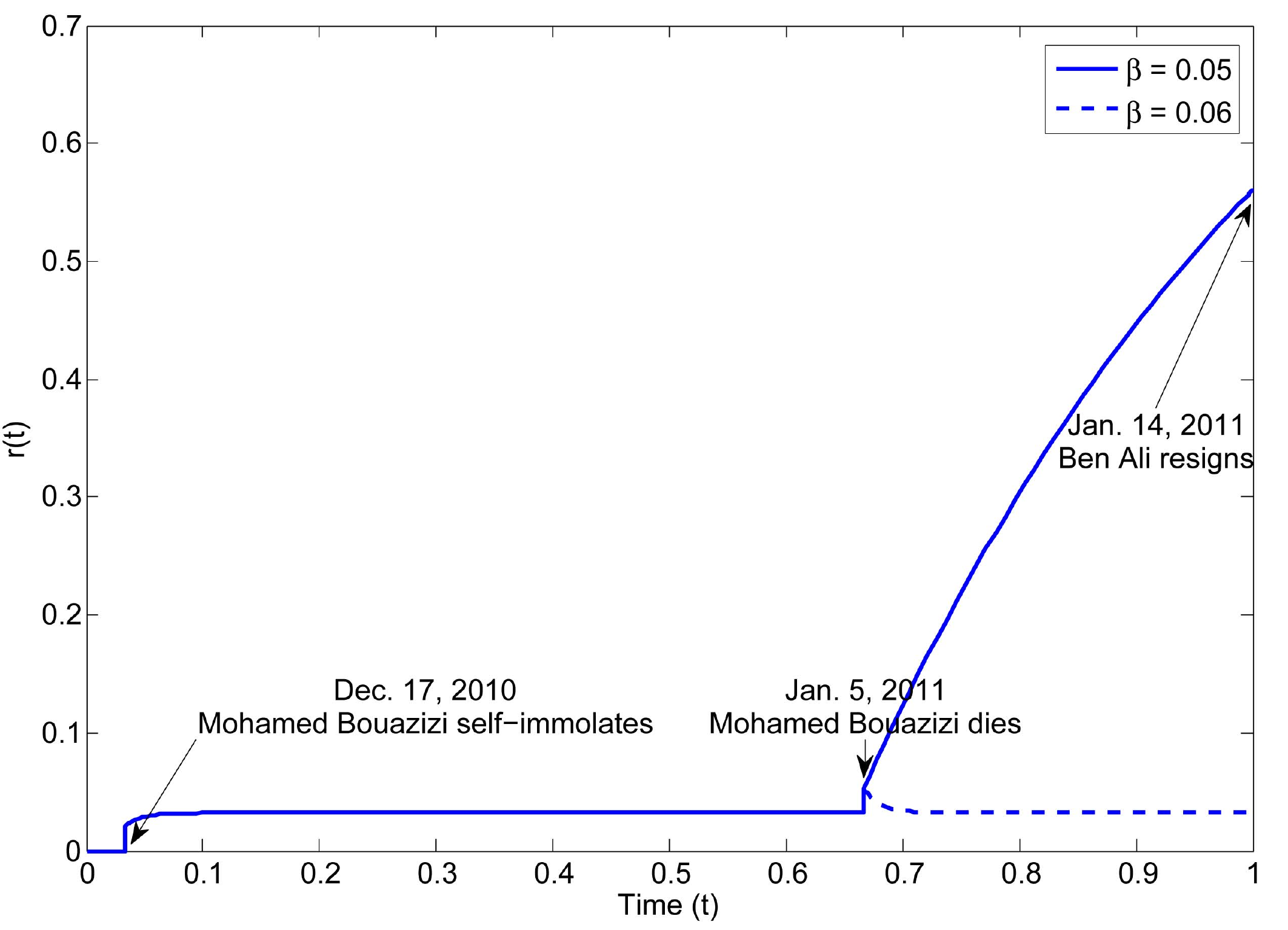}
		\caption{Meta-stable police state (region IIIe) with $\alpha=0.98$}
		\label{fig:TunisiaTimeAlphaB}
	\end{subfigure}
	\caption{The effect of increasing $\alpha$ on the behaviour of solutions to \eqref{eq:ODE} with $\beta\in\{0.05,0.06\}$, $c_1=2.30$, $c_2=30\log(10)$, and subject to shocks $\Delta r_1 = \Delta r_2=0.021$ occurring at $t = \frac{1}{30}$ (December 17, 2010) and $t=\frac{20}{30}$ (January 5, 2011).  Panel (b) with $\beta=0.05$ corresponds qualitatively to Figure 2 of \citet{HowardEtAl11}.}
	\label{fig:TunisiaTimeAlpha}
\end{figure}

\begin{figure}[h]
	\centering
	\begin{subfigure}[h]{3in}
		\centering
		\includegraphics[width=3in]{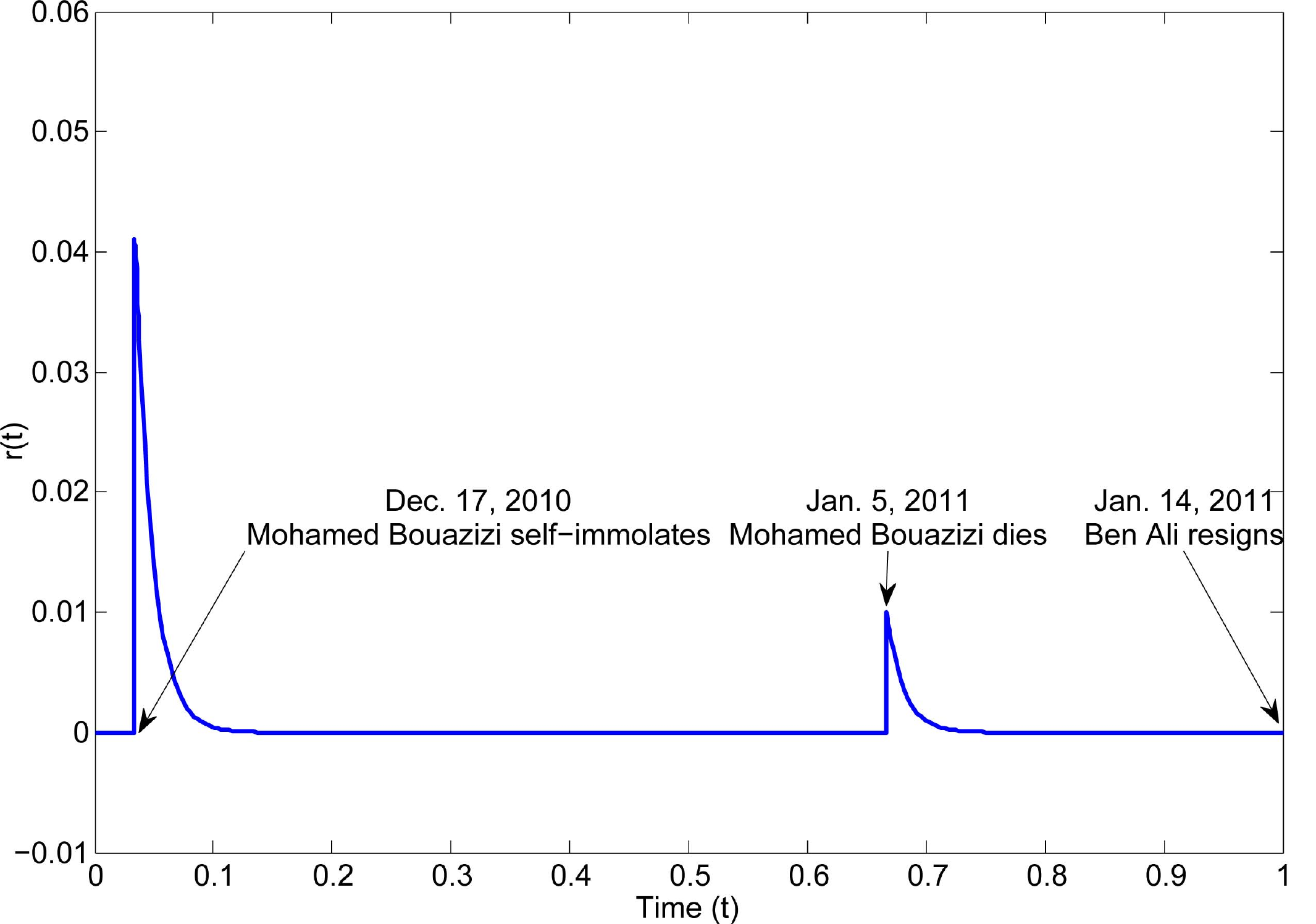}
		\caption{Stable police state (III0) with $c_1=2.30$.}
		\label{fig:TunisiaTimeC1A}
	\end{subfigure}
	\hspace{5mm}
	\begin{subfigure}[h]{3in}
		\centering
		\includegraphics[width=3in]{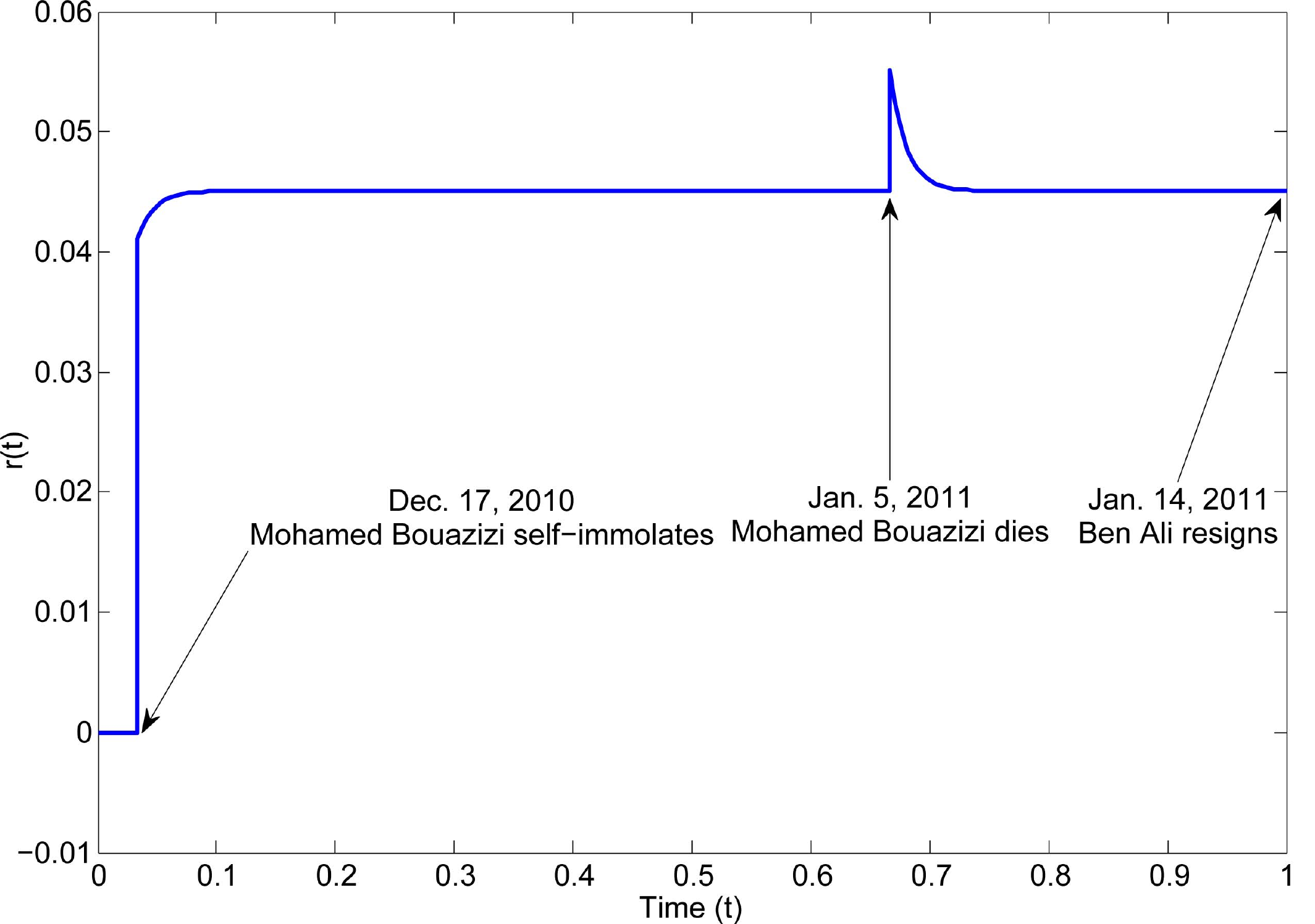}
		\caption{Meta-stable police state (IIIe) with $c_1=3.26$.}
		\label{fig:TunisiaTimeC1B}
	\end{subfigure}\\
	\begin{subfigure}[h]{3in}
		\centering
		\includegraphics[width=3in]{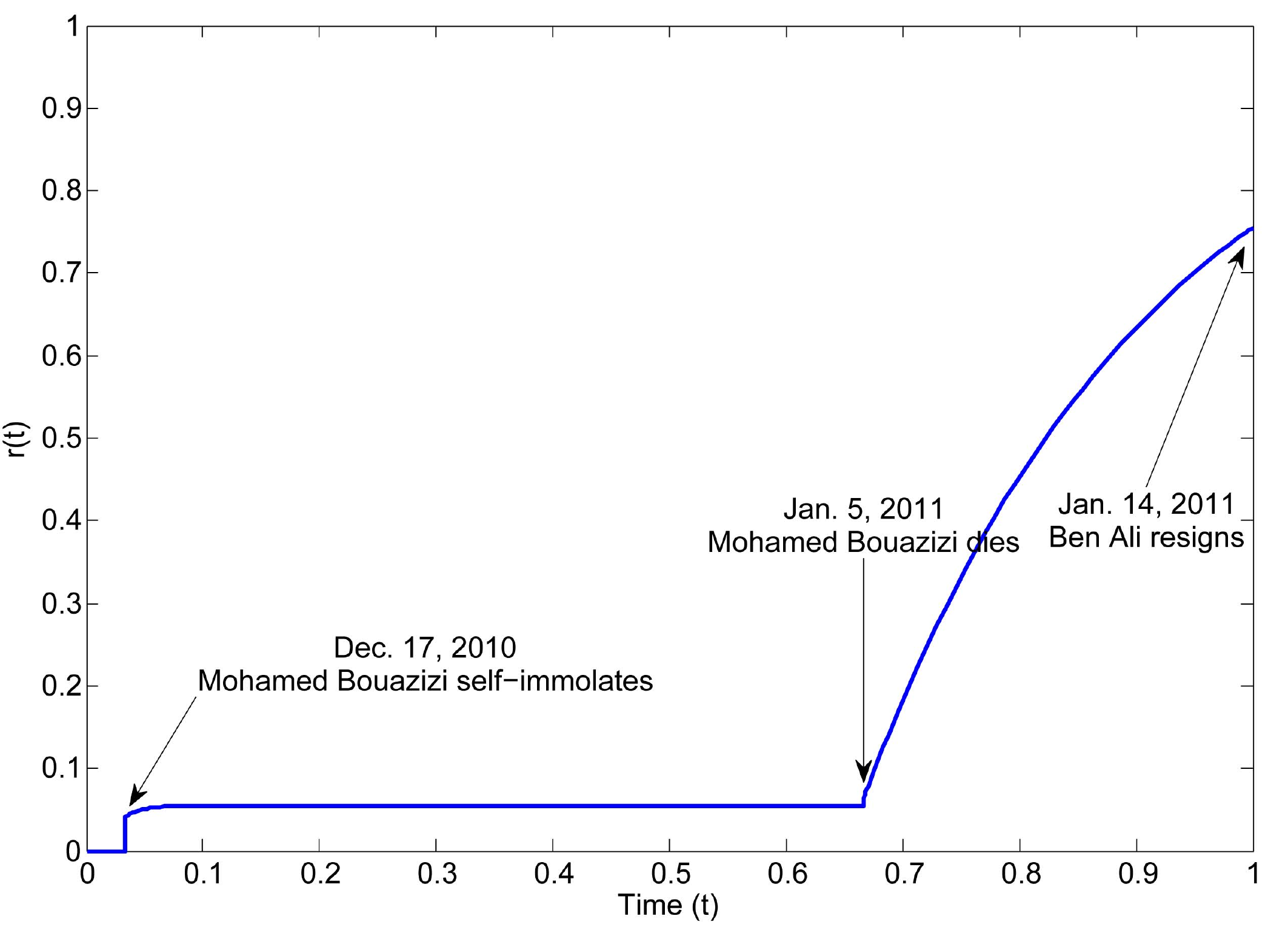}
		\caption{Meta-stable police state (IIIe) with $c_1=4.02$.}
		\label{fig:TunisiaTimeC1C}
	\end{subfigure}
	\hspace{5mm}
	\begin{subfigure}[h]{3in}
		\centering
		\includegraphics[width=3in]{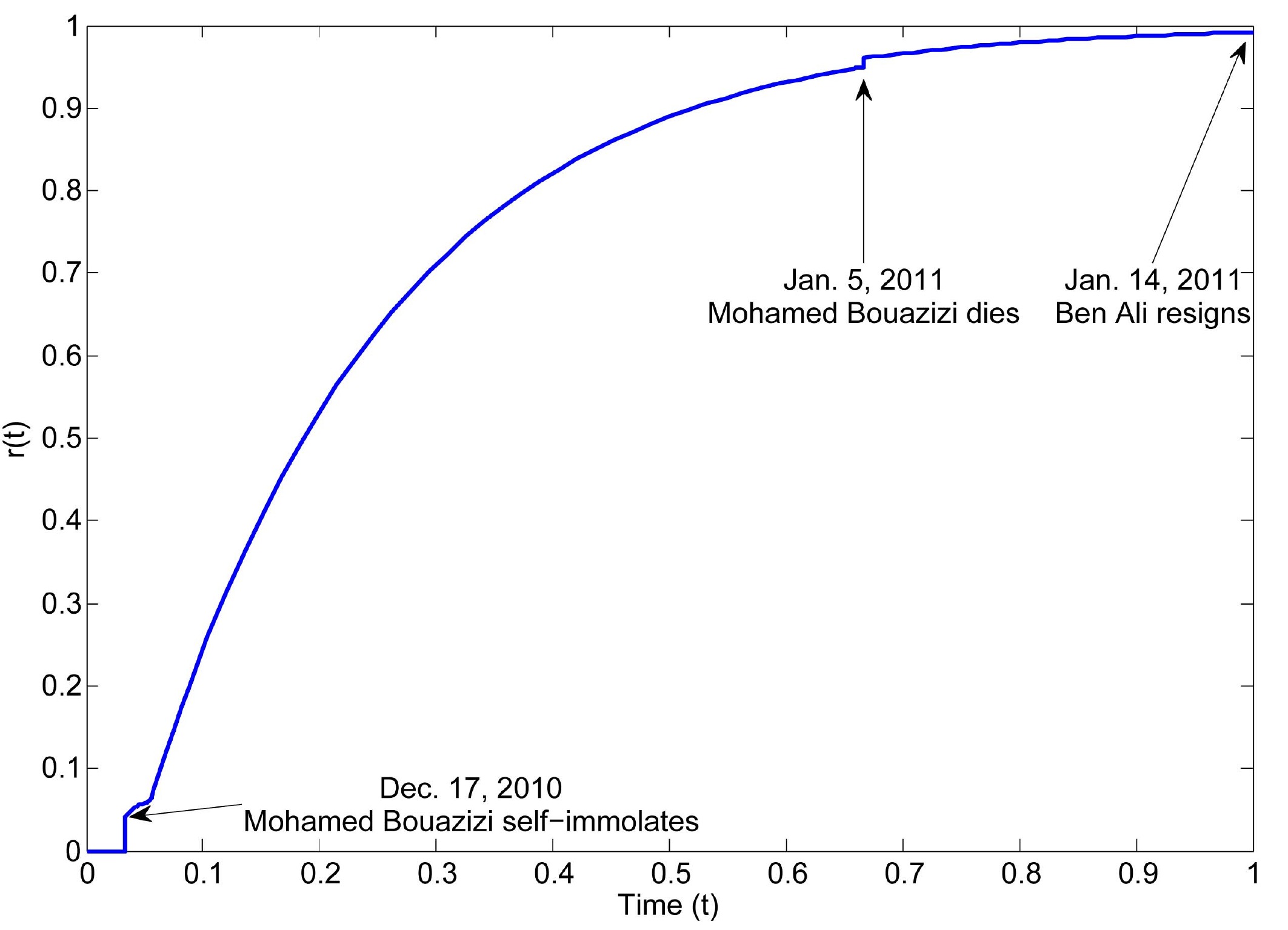}
		\caption{Unstable police state (III1) with $c_1=4.80$}
		\label{fig:TunisiaTimeC1D}
	\end{subfigure}
	\caption{The effect of increasing $c_1$ on the behaviour of solutions to \eqref{eq:ODE} with $\alpha=0.96$, $\beta=0.06$, $c_2=30\log(10)$, and subject to shocks $\Delta r_3=0.041$ and $\Delta r_4=0.01$ occurring at $t = \frac{1}{30}$ (December 17, 2010) and $t=\frac{20}{30}$ (January 5, 2011).  Panel (c) corresponds qualitatively to Figure 2 of \citet{HowardEtAl11}.}
	\label{fig:TunisiaTimeC1}
\end{figure}

\paragraph{}  Another way to visualize the information presented in Figure \ref{fig:TunisiaLine} is to plot time traces of solutions to \eqref{eq:ODE} for different choices of parameters.  The perturbations mentioned above, which are required to move from the basin of attraction of one equilibrium to another, are delivered by shocks to the fraction of revolutionaries ($\Delta r$).  Specifically, although one shock of sufficient magnitude would perturb a solution from total state control ($r=0$) to the basin of attraction of the realized revolution ($r=1$), we choose to consider a scenario with two smaller shocks.  This is because the probability of a shock likely decreases rapidly with the its magnitude, and therefore, we are significantly more likely to encounter two smaller shocks as opposed to one large shock.  In fact, the two shocks are chosen to coincide with the self-immolation (December 17, 2010) and death (January 5, 2011) of Mohamed Bouazizi which, respectively, set in motion the Tunisian revolution and immediately preceded a spike in online conversations about freedom and revolution \citep{HowardEtAl11} in the lead-up to the resignation of President Ben Ali (January 14, 2011).  Using this time scale, Ben Ali resigns at $t=\frac{29}{30}$.  In keeping with the calculations done at the beginning of Section \ref{sec:Interp}, we choose $c_2=30\log(10)\approx69.1$.  Figure \ref{fig:TunisiaTimeAlpha} shows the effect of increasing $\alpha$ by solving \eqref{eq:ODE} with 

$$\alpha=\{0.96,0.98\},\mbox{ } \beta\in\{0.05,0.06\}, c_1=2.30 \mbox{ }(c^*\approx 0.0322), \mbox{ }r(0)=0,$$

and subject to shocks $\Delta r_1 = 0.021$ and $\Delta r_2=0.021$ occurring on December 17, 2010, ($t=\frac{1}{30}$) and January 5, 2011, ($t=\frac{20}{30}$), respectively.  Figure \ref{fig:TunisiaTimeC1} demonstrates the effect of increasing $c_1$ by solving \eqref{eq:ODE} with 

$$\alpha=0.96,\mbox{ } \beta=0.06, \mbox{ }c_1\in\{2.30, 3.26, 4.02, 4.80\}\mbox{ } (c^*\in\{0.0322,0.0451,0.0550, 0.0650\}),$$ 

and subject to shocks $\Delta r_3=0.041$ and $\Delta r_4=0.01$ occurring on December 17, 2010, ($t=\frac{1}{30}$) and January 5, 2011, ($t=\frac{20}{30}$), respectively. 

\paragraph{}  As with Figure \ref{fig:TunisiaAlpha}, Figure \ref{fig:TunisiaTimeAlpha} illustrates how the dynamics of the model change as $\alpha$ increases.  Because solutions shown in Figure \ref{fig:TunisiaTimeAlphaA} have parameters in the stable police state (region III0), leaving the basin of attraction of total state control ($r=0$) requires a shock $\Delta r>\beta$.  Thus, the solutions in Figure \ref{fig:TunisiaTimeAlphaA} decaying to total state control is consistent with $\Delta r_1=\Delta r_2 = 0.021<\beta$ for both $\beta=0.05$ and $\beta = 0.06$.  Increasing the visibility from $\alpha=0.96$ to $\alpha = 0.98$ decreases the basin of attraction of total state control and yields Figure \ref{fig:TunisiaTimeAlphaB}.  In this case, the shock $\Delta r_1$ is sufficient to move solutions from total state control to the basin of attraction of civil unrest ($r=c^*$).  However, as we saw in Figure \ref{fig:TunisiaAlpha}, for a second shock to perturb the solution from civil unrest to the basin of attraction of the realized revolution ($r=1$) requires $\Delta r> \beta- c^*$.  This is the case for $\beta = 0.05$ but not $\beta = 0.06$.  

\paragraph{}  Similar to Figure \ref{fig:TunisiaLine}, Figure \ref{fig:TunisiaTimeC1} illustrates how the dynamics of the model change as $c_1$ increases.  Figure \ref{fig:TunisiaTimeC1A} shares the same parameters as Figure \ref{fig:TunisiaTimeAlphaA} and shows solutions in the stable police state (region III0) for which the two shocks, now taken with strength $\Delta r_3 = 0.041$ and $\Delta r_4=0.01$, are insufficient to leave the basin of attraction of total state control ($r=0$).  As the enthusiasm of protesters increases from $c_1=2.30$ ($c^*=0.0322$) to $c_1=3.26$ ($c^*=0.0451$) we go from a stable police state to a meta-stable police state (region IIIe). This situation is illustrated in Figure \ref{fig:TunisiaTimeC1B}, where $\Delta r_3 = 0.041$ is now sufficient to move the solution from total state control to the basin of attraction of civil unrest ($r=c^*$), however, $\Delta r_4$ is still too small relative to $\beta-c^*$ to perturb the solution from civil unrest to the basin of attraction of the realized revolution ($r=1$), i.e. $\Delta r_4 < \beta - c^*$.  Continuing to increase the enthusiasm parameter to $c_1 = 4.02$ ($c^* = 0.0550$) maintains parameters in the meta-stable police state.  Here the magnitude of a shock needed to leave total state control for the basin of attraction of civil unrest remains unchanged, whereas the magnitude of a shock needed to perturb a solution from civil unrest to the basin of attraction of the realized revolution is decreased.  This is illustrated in Figure \ref{fig:TunisiaTimeC1C}, where $\Delta r_4$ is now sufficient to propel the solution from civil unrest to the basin of attraction of the realized revolution.  Finally, Figure \ref{fig:TunisiaTimeC1D} shows the situation where further increased enthusiasm, $c_1=4.80$ ($c^*=0.0650$), has moved parameters from a meta-stable police state to an unstable police state (region III1).  Observe that for the unstable police state the shock $\Delta r_3$ is sufficient on its own to move the solution from total state control to the basin of attraction of the realized revolution.

\paragraph{}  We summarize findings from Figures \ref{fig:TunisiaLine}-\ref{fig:TunisiaTimeC1} by observing that (a) an increase in visibility ($\alpha$) decreases the size of shock needed to leave total state control ($r=0$) for the basin of attraction of civil unrest ($r=c^*$), and (b) an increase in enthusiasm ($c_1$) decreases the size of the shock needed to leave civil unrest for the basin of attraction of the realized revolution $(r=1)$.  Indeed, for sufficiently high $c_1$ ($c^*$) no second shock is needed at all.  When considering a revolution, however, we must consider effects of both visibility and enthusiasm simultaneously, because as discussed above we expect the Internet and social media to increase both $\alpha$ and $c_1$.

\paragraph{}  When comparing the plots in Figures \ref{fig:TunisiaTimeAlphaB} and \ref{fig:TunisiaTimeC1C} to data presented in Figure 2 of \citet{HowardEtAl11}\footnote{This figure is captioned ``Percent of Tunisian Blogs With Posts on Politics, By Keyword''.} we find a qualitative match.  Howard's Figure 2 shows an initial increase in posts related to the economy following Mohamed Bouazizi's self immolation starting on December 17, 2010, which is sustained until Mohamed Bouazizi's death on January 5, 2011, when we see the start of a second increase in posts related to Ben Ali.  This second increase peaks on January 14, 2011, when Ben Ali resigns.  \citet{HowardEtAl11} observes that street protests continued well after the departure of Ben Ali until at least the resignation of Prime Minister Mohamed Ghannouchi on February 27, 2011.  This appears consistent with the behaviour in Figures \ref{fig:TunisiaTimeAlphaB} and \ref{fig:TunisiaTimeC1C} where we see an initial spike leading to sustained engagement and a second spike leading to full-on revolution that continues to gain momentum well after the January 14, 2011, ($t=1$) resignation of Ben Ali.

\paragraph{}  Now that we have a better understanding of how adoption of the Internet and social media may factor into our model, let us consider aspects of the Jasmine revolution that we can investigate using our model by addressing the following two questions.

\begin{enumerate}
	\item  How can a small number of active social media users and relatively low Internet penetration have a dramatic effect on the stability of a regime?, and
	\item  How is it that regimes manage to seem so stable until the revolution is underway?
\end{enumerate}

\paragraph{}\emph{Question 1:}  To explain how only a small number of social media users can have a significant impact on the likelihood of a full-blown revolution, consider the effect of the Internet and social media usage solely on visibility ($\alpha$).  In the example of Figure \ref{fig:TunisiaTimeAlpha} a small increase from $\alpha=0.96$ to $\alpha=0.98$ reduces the size of shock in $r$ necessary to leave total state control ($r=0$) from $\Delta r = 0.04$ to $\Delta r = 0.02$.  If shocks occur distributed according to some probability distribution, then it is reasonable to assume that shocks of sufficient magnitude to mobilize large fractions of the population lie in the tail of this distribution.  For many reasonable probability distributions satisfying this criterion, halving the size of shock necessary to trigger a revolution more (and potentially much more) than doubles the likelihood of a revolution occurring in any given amount of time.  Compounding this phenomenon is how, when $\alpha$ initially increases beyond $1 - c^*$, the basin of attraction of total state control ($r=0$) shrinks from $(0,\beta)$ to $(0,1-\alpha]$ in a discontinuous fashion.  So, a small increase in $\alpha$ can have a very large impact on the expected amount of time one has to wait until a revolution is triggered.

\paragraph{}\emph{Question 2:}  Increasing either enthusiasm ($c_1$) or visibility ($\alpha$) eventually decreases the size of the basin of attraction of total state control ($r=0$).  This undermines the regime by decreasing the size  of shock necessary to trigger a revolution.  However, since (a) $r=0$ always remains a locally asymptotically stable equilibrium, (b) large shocks are exceedingly rare, and  (c) determining the exact values of the parameters in a model like ours is very difficult, the exact size of shock necessary to trigger a revolution is impossible to determine until such a shock occurs.  It follows that for someone observing a regime before and after the adoption of social media there would be few, if any, outward signs of instability: the regime appears stable until it isn't.

\subsection{Case Study: Egypt}
\label{sec:Egypt}
\paragraph{}  In applying our model to the Tunisian case study we included no factors external to the model, other than shocks to the number of protesters, $r$.  This was sufficient to produce a qualitative match between the development of the Tunisian revolution and the model.  In the case of the Egyptian revolution, however, there are certain singular events not captured by our model that would have a significant effect on the parameters of the model.  For example,

\begin{enumerate}
	\renewcommand{\labelenumi}{(\alph{enumi})}
	\item news of a successful revolution in Tunisia likely raised the enthusiasm ($c_1$) and visibility ($\alpha$) parameters in Egypt by causing Egyptians to discuss and re-evaluate the strengths (and weaknesses) of the regime, the discontent of the general population, and their chances of success \citep{BBC11, Pollock11, ZhuoEtAl11},
	\item the intervention of the Egyptian military on behalf of protesters in Tahrir Square \citep{BBC11, Said12} likely lowered policing capacity ($\beta$) and efficiency ($c_2$) directly by immediately curtailing the regime's policing capacity, and
	\item the January 28 - February 1, 2011 Internet disruptions \citep{Dunn11, HowardEtAl11, BBC11} lowered the visibility ($\alpha$) temporarily.
\end{enumerate}

\paragraph{}  This suggests to consider changes in the model parameters during the course of the revolution.  First we establish the following rough timeline \citep{BlightEtAl12, Dunn11, BBC11, HowardEtAl11, Pollock11, Said12}.

\begin{itemize}
	\item December 17, 2010: Mohamed Bouazizi self-immolates in Tunisia.
	\item January 14 - 15, 2011: Ben Ali flees Tunisia and an interim government is established.
	\item January 25, 2011: Day of Protest in Tahrir Square, Egypt.
	\item January 26, 2011: Police clear Tahrir Square.
	\item January 28, 2011: Protesters occupy Tahrir Square, Mubarak addresses nation, major Internet disruptions begin.
	\item February 1, 2011: President Obama withdraws support for Mubarak regime, army refuses to act against protesters, major Internet disruptions end.
	\item February 2, 2011: State vandals and thugs attack protesters in Tahrir Square, army officers intervene on behalf of protesters.
	\item February 11, 2011: Mubarak resigns.
\end{itemize}

\begin{figure}[h!]
	\centering
	\includegraphics[width=6in]{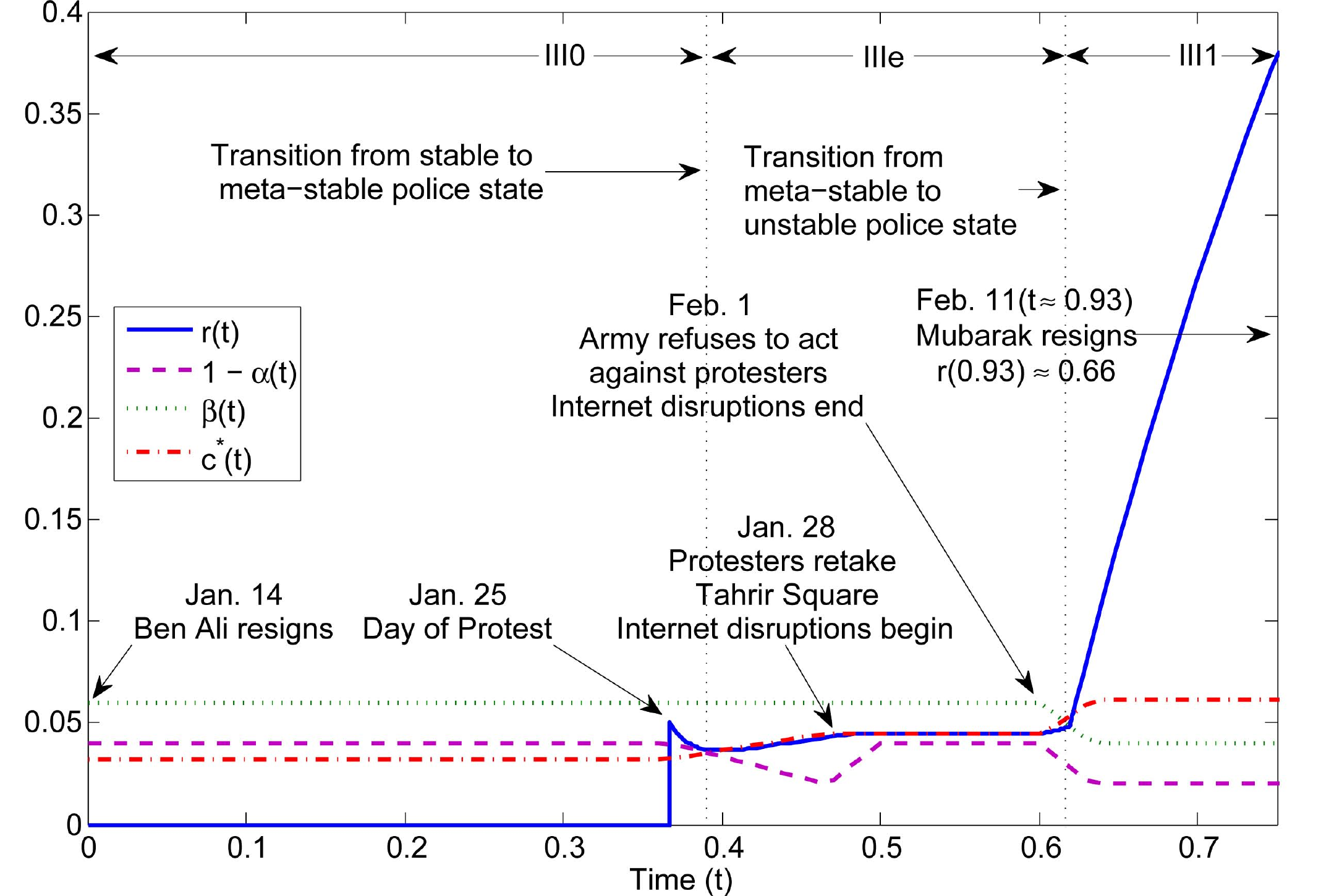}
	\caption{A possible scenario modelling the Egyptian revolution using dynamic parameters and an initial shock.  The model solution corresponds qualitatively to Figure 4 of \citet{HowardEtAl11}.}
	\label{fig:Egypt}
\end{figure}

\paragraph{}  We now outline a possible interpretation for the revolution in Egypt according to our model.  The events of January 14-15 likely increased both $\alpha$ and $c_1$ but not by enough to move Egypt from a stable (region III0) to a meta-stable (region IIIe) police state.  Because the initial shock of the Day of Protest on January 25 is of insufficient magnitude to perturb the system from total state control ($r=0$) to the basin of attraction of the realized revolution ($r=1$), i.e. the magnitude of the shock is less than $\beta$, Egyptian police are able to clear Tahrir Square.  Nevertheless, this initial protest is sufficiently large to generate considerable coverage on both social and traditional media \citep{Alterman11, HowardEtAl11, ZhuoEtAl11}.  An increased awareness in the general population of the current level of dissatisfaction with the regime through increased consultation of Internet and satellite TV sources then resulted in a further increase of $\alpha$.  In addition, protesters' initial success leads to increased enthusiasm and experience which may increase $c_1$.  At some point between January 25 and 28 Egypt becomes a meta-stable police state and the revolution converges to the equilibrium of civil unrest ($r=c^*$).  Meanwhile, Internet disruptions between January 28 and February 1 temporarily suppress $\alpha$.  Finally, the decisions made by the army in favor of the protesters on February 1-2 lower both $\beta$ and $c_2$.  Egypt then becomes an unstable police state (region III1) and the revolution proceeds to completion without the need of a second shock.  We illustrate this scenario in Figure \ref{fig:Egypt} where $t=0$ is taken to be January 14, $r(0)=0$, parameters\footnote{Parameters are chosen so that they conform with the scenario presented above, but their precise values are chosen rather arbitrarily within these constraints for this proof-of-concept scenario analysis.  We do not make any attempt to model the time-dependence of the parameters explicitly.} are given in Table \ref{tab:EgyptParam} and in equations \eqref{eq:EgyptParamA}-\eqref{eq:EgyptParamD}, and the initial shock occurs on January 25 ($t=\frac{11}{30}$) with $\Delta r = 0.05$.

\begin{table}[h]
	\centering
	\begin{tabular}{llllc}
		Date & Time (t) &  Event & Parameters Affected & Region\\
		\hline
		Jan. 14-15 & $0$ to $\frac{1}{30}$ & $\cdot$ Ben Ali resigns and an & $\cdot$ Set $\alpha$ and $c_1$ to & III0\\
		&& \hspace{2mm} interim government is & \hspace{2mm} $\alpha = 0.96$ and $c_1=2.30$\\
		&& \hspace{2mm} formed in Tunisia& $\cdot$ Initially $\beta = 0.06$ and\\
		&&& \hspace{2mm}  $c_2 = 69.1$\\
		Jan. 25 & $\frac{11}{30}$ & $\cdot$ Day of Protest in Tahrir Square & $\cdot$ Shock to $r$ of $\Delta r= 0.05$ & III0\\
		Jan. 25-28 & $\frac{11}{40}$ to $\frac{14}{30}$ & $\cdot$ Lead up to Jan. 28 protest & $\cdot$ Increase $\alpha$ and $c_1$ & III0 - IIIe\\
		&&& \hspace{2mm} (linearly) to $\alpha = 0.98$\\
		&&& \hspace{2mm} and $c_1 = 3.26$\\
		Jan. 28 - Feb. 2 & $\frac{14}{40}$ to $\frac{19}{30}$ & $\cdot$ Temporary Internet & $\cdot$ Suppress\footnotemark $\alpha$ from & IIIe \\
		&& \hspace{2mm} disruptions & \hspace{2mm} $\alpha=0.98$ to $\alpha=0.96$\\
		Feb. 1-2 & $\frac{18}{30}$ to $\frac{19}{30}$ & $\cdot$ Army intervenes on behalf & $\cdot$ Decrease\footnotemark[7] $\beta$ and $c_2$ to & IIIe - III1\\
		&&\hspace{2mm} of protesters & \hspace{2mm} $\beta = 0.04$ and $c_2 = 50.0$
	\end{tabular}
	\caption{Significant events during the Egyptian revolution and their impact on parameters $\alpha$, $\beta$, $c_1$, and $c_2$.  The parameter profiles as a function of time are described in the text and are plotted in Figure \ref{fig:Egypt}.}
	\label{tab:EgyptParam}
\end{table}

\footnotetext[7]{Transitions are assumed to occur linearly over the period of one day.}

\newpage

\paragraph{}  The parameter profiles are given as a function of time by

\begin{subequations}
	\label{eq:EgyptParam}
	\begin{align}
		\alpha(t) &= \left\{ \begin{array}{lrl} 0.96 & 0\leq t<\frac{11}{30} & \mbox{Prior to Jan. 25}\\ 
			0.96\frac{14-30t}{3} + 0.98\frac{30t-11}{3} & \frac{11}{30} \leq t < \frac{14}{30} & \mbox{Jan. 25-28}\\
			0.98(15-30t) + 0.96(30t-14) & \frac{14}{30} \leq t <\frac{15}{30} & \mbox{Jan. 28-29}\\
			0.96 & \frac{15}{30}\leq t<\frac{18}{30} & \mbox{Jan. 29 - Feb. 1}\\
			0.96(19-30t) + 0.98(30t-18) & \frac{18}{30}\leq t<\frac{19}{30} & \mbox{Feb. 1-2}\\
			0.98 & t\geq \frac{19}{30} & \mbox{Feb. 2 onwards} \end{array} \right.,
			\label{eq:EgyptParamA}\\
		c_1(t) &= \left\{ \begin{array}{lrl} 2.30 & 0\leq t<\frac{11}{30}& \mbox{Prior to Jan. 25}\\ 
			2.30\frac{14-30t}{3} + 3.26\frac{30t-11}{3} & \frac{11}{30} \leq t < \frac{14}{30}& \mbox{Jan. 25-28}\\
			3.26 & t\geq \frac{14}{30} & \mbox{Jan. 28 onwards}  \end{array} \right.,\\	
		\beta(t) &= \left\{ \begin{array}{lrl} 0.06 & 0\leq t<\frac{18}{30} & \mbox{Prior to Feb. 1}\\ 
			0.06(19-30t) + 0.04(30t-18) & \frac{18}{30} \leq t < \frac{19}{30} & \mbox{Feb. 1 - 2}\\
			0.04 & t\geq \frac{19}{30} & \mbox{Feb. 2 onwards} \end{array} \right., \mbox{ and}\\
		c_2(t) &= \left\{ \begin{array}{lrl} 69.1 & 0\leq t<\frac{18}{30} & \mbox{Prior to Feb. 1}\\ 
			 69.1(19-30t) + 50.0(30t-18) & \frac{18}{30} \leq t < \frac{19}{30} & \mbox{Feb. 1 - 2}\\
			50.0 & t\geq \frac{19}{30} & \mbox{Feb. 2 onwards}  \end{array} \right..
			\label{eq:EgyptParamD}
	\end{align}
\end{subequations}

\paragraph{}  As in the previous case study of Tunisia, the model scenario illustrated by Figure \ref{fig:Egypt} agrees qualitatively with observed data, in particular the data presented in Figure 4 of \citet{HowardEtAl11}\footnote{This figure is captioned ``Logged Number of Tweets on \#egypt, by Location''.}, which shows rapid growth in Tweets related to the Egyptian revolution culminating in a spike of Twitter activity on January 25.  After this spike, activity declines and then stabilizes between January 25 and January 28.  Finally, after restrictions on the Internet are lifted on February 1 there is rapid growth in Twitter activity until Mubarak resigns on February 11.  Figure \ref{fig:Egypt} shows that our model is again capable of generating behaviour observed in \citet{HowardEtAl11}.

\paragraph{}  We take this opportunity to address the third question posed in Section \ref{sec:Motive},

\begin{enumerate}
	\setcounter{enumi}{2}
	\item  Why did the January 28 - February 1 Internet shutdown in Egypt not have a greater inhibitory effect on protests?
\end{enumerate}

\paragraph{}\emph{Question 3:}  In our model, attempting to suppress a revolution by reducing its visibility ($\alpha$) will only be successful if the regime can make the visibility threshold ($1-\alpha$) greater than the current size of the revolution.  Figure \ref{fig:Egypt} illustrates the scenario where the regime fails in its attempt to suppress a revolution because it does not manage to sufficiently decrease visibility of protesters.  Thus, the January 28 - February 1 Internet shutdown had little effect on the ongoing protests because the protests had become sufficiently large that individuals no longer needed the Internet to be aware of them.  We make the additional observation that, with respect to the example of Figure \ref{fig:Egypt}, increasing $\beta$ might prevent the transition from a meta-stable police state (region IIIe) to an unstable police state (region III1), but would not prevent the transition from a stable police state (region III0) to a meta-stable police state.  As such, according to our model, the regime is not able to regain total state control ($r=0$) by increasing its policing capacity, $\beta$.  The only other option open to the regime to suppress the revolution would be to decrease $c^*$ by increasing its policing efficiency ($c_2$).  Increasing policing efficiency rapidly is difficult, however, because improving training of security forces, intelligence, investments in infrastructure and crowd control, etc... take a significant amount of time.  Moreover, as we discussed in Section \ref{sec:Tunis} attempting to increase policing efficiency by use of police brutality risks being exposed by new media, thus inducing otherwise apolitical individuals to join the revolution \citep{BBC11, Schneider11}.

\paragraph{}  Sections \ref{sec:Tunis} and \ref{sec:Egypt} have presented ways in which our model can be applied to the revolutions of Tunisia and Egypt.  These two case studies have shown how our model, taking into account only the visibility of protesters and policing capacity of the regime, can produce results that are qualitatively consistent with observations and data presented in \citet{HowardEtAl11}.  This indicates that the visibility of protesters and policing capacity of the regime may be essential factors for the dynamics of the Arab Spring revolutions.  Since the visibility of protesters is primarily affected by media and communications technologies, it follows that through the case studies presented above our model supports the view that the inernet and social media played a critical role in the Arab Spring 2010-2011 \citep{HowardEtAl11, LotanEtAl11, Alterman11, BBC11, KhamisVaughn11, Pollock11, Saletan11, Shirky11, Stepanova11, ZhuoEtAl11}.

\subsection{Case Studies: Iran, China, and Somalia}
\label{sec:Other}
\paragraph{}  Above we posed the question

\begin{enumerate}
	\setcounter{enumi}{3}
	\item  Why is it that some regimes fall in a matter of weeks, others fight to a stalemate, and still others survive relatively unscathed?
\end{enumerate}

In our previous discussions we have seen that when a regime is subject to a shock the outcome depends on the balance between the four parameters in our model: the visibility ($\alpha$) and enthusiasm ($c_1$) of protesters, and the capacity ($\beta$) and efficiency ($c_2$) of the police.  We further explore this point by considering the cases of Iran, China, and Somalia.

\paragraph{}\emph{Question 4:}  The protests following Iran's 2009 election, dubbed the ``Green Revolution'', were ultimately put down by the regime despite widespread use of social media technology.  In particular, \citet{BurnsEltham09} emphasize the response of Iran's Revolutionary Guard and paramilitary force, the Basij, which unleashed a brutal crackdown in part by using Twitter to ``hunt down and target Iranian pro-democracy activists''.  Unfortunately, this result was partly caused by Haystack, a poorly vetted anti-censorship software promoted by the US government, which 

\begin{quote}
	``not only failed at its goal of hiding messages from governments but also made it, in the words of one analyst, `possible for an adversary to specifically pinpoint individual users.' '' \citep{Shirky11}
\end{quote}

The large amount of resources that were available to the Iranian regime is consistent with a large value of $\beta$.  The ability of the regime to harness social media to suppress protests, and to mobilize a well equipped and motivated security force is consistent with a large $c_2$.  Moreover, since at the time of the Green Revolution social media was still in its infancy\footnote{Facebook was launched in 2004 and was still an invitation-only service in 2005\citep{Phillips07}.  Youtube was founded in early 2005 \citep{Youtube} and Twitter was not founded until the spring of 2006 \citep{Picard11}.}, $\alpha$ and $c_1$ were unlikely to have felt the full impact of these new communications technologies.  Contributing to the lower values of $\alpha$ and $c_1$ relative to the two case studies presented above, is the lack of both experience with and prior examples of protests.  When $c_1$ is small and $c_2$ is large, $c^*$ is small.  So, in our model a 2009-era Iran with small $c^*$ and $\alpha$, and large $\beta$, specifically $c^*<1-\alpha<\beta$, is a stable police state (region III0), which is consistent with the failure of the Green Revolution.  In the future, more sophisticated methods for evading government detection and identification may increase $c_1$ at the expense of $c_2$.  Continued growth in Internet availability will increase $\alpha$ by increasing visibility, while economic sanctions may increase $\alpha$ by stoking dissatisfaction with the regime.  Examples of successful revolutions in the Arab world and contact/support from successful revolutionaries from abroad will also increase $\alpha$ and $c_1$.  If these factors manage to out pace the evolution of the regime's police forces (including their technological abilities) then it can be expected that Iran may pass into the meta-stable (region IIIe) or unstable (region III1) police state regions in the future.  A subsequent shock, perhaps due to another highly contested election, may then trigger a revolution.

\begin{figure}[h!]
	\centering
	\includegraphics[width=3in]{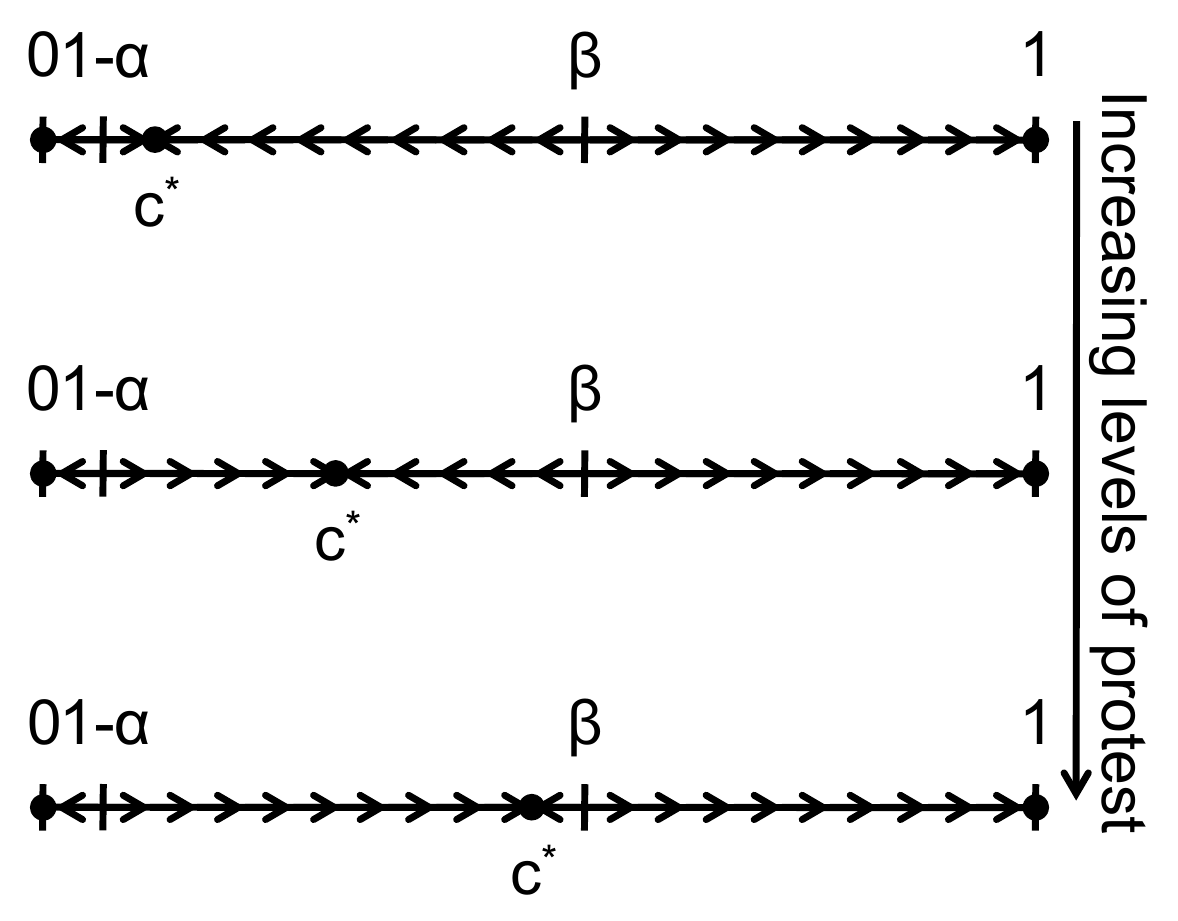}
	\caption{Effect of increasing $c^*$ when $\alpha+\beta>1$ and $\frac{c_1}{c_1+c_2}=c^*>1-\alpha$ (region IIIe).  Increasing $c_1$ causes $c^*$ to increase as well.  For small $c_1$ ($c^*$) we see small numbers of protests that do not endanger the regime.  Growing $c_1$ ($c^*$) causes the size/number of protests to increase.  Eventually, for large $c_1$ ($c^*$) the regime becomes endangered as the probability of a shock sufficient to trigger a full-scale revolution becomes significant.}
	\label{fig:China}
\end{figure}

\paragraph{}  While the current regime in China differs from the pre-revolutionary regimes in Tunisia and Egypt in many aspects, it is interesting to consider how our model may apply to China in terms of the influence of state control on the media and the Internet, and police control of dissident opinion.  The number of ``mass group incidents'' reported annually in China has been rising consistently for at least two decades \citep{Wedeman09}.  Being constantly subject to low but rising levels of protest may correspond to the civil unrest equilibrium ($r=c^*$) in region IIIe of our model, which we have called the meta-stable police state region.  In our model rising levels of protest would correspond to rising $c^*$, where in this particular case, an increasing $c^*$ would seem to be the result of an increase in the enthusiasm ($c_1$) of protesters and not a decrease in the efficiency of the regime ($c_2$), except perhaps in terms of Internet censorship.  Previously, we argued that regimes in the meta-stable police state were potentially at an increased risk of revolution depending on the balance of $\alpha$, $c^*$, and $\beta$.  Specifically, we stated that increasing $c_1$, and hence $c^*$, decreases (increases) the magnitude of shock needed to go from $r=c^*$ to the basin of attraction of $r=1$ ($r=0$), see Figure \ref{fig:China}.  A continued rise of $c^*$ (via increasing $c_1$) and a sudden or systematic rise in $\alpha$ through Internet and social media exposure in China may eventually result in increasing the chance of a successful revolution.  How soon this would occur is not possible to say, since we are unable to determine accurately the magnitudes of $1-\alpha$, $c^*$, and $\beta$.

\paragraph{}  Finally, we consider the case of Somalia, which is widely considered to have been a failed state for more than 20 years.  During this period the country has seen the rise and fall of many local authorities and attempts at re-establishing a national government \citep{Menkhaus07}.  The failed state region (region II) features low $\alpha$ (weak media) and low $\beta$ (weak government).  Low visibility ($\alpha$) prevents individuals from joining any popular movements and low policing capacity ($\beta$) prevents the government from reigning in existing movements.  This results mathematically in an uncountable number of equilibria contained in a subinterval of $[0,1]$, and is consistent with the slow and erratic rise and fall of local militia and a succession of weak central governments.  Our model predicts that a successful national state could arise from either (a) improving policing capacity of the transitional government (increasing $\beta$), or (b) increasing social cohesion and the capacity of the media in Somalia (increasing $\alpha$).  Interestingly, due to the lack of central authority, as well as permanent infrastructure, Somalia has developed a sophisticated and affordable telecommunications sector \citep{Feldman07}, which may mean that an increased $\alpha$ is not unrealistic.

\section{Discussion and Conclusion}
\label{sec:Conc}
\paragraph{}  We began this paper by asking four questions related to the revolutions of the Arab Spring.

\begin{enumerate}
	\item  How can a small number of active social media users and relatively low Internet penetration have a dramatic effect on the stability of a regime?,
	\item  How is it that regimes manages to seem so stable until the revolution is underway?,
	\item  Why did the January 28 - February 1 Internet shutdown in Egypt not have a greater inhibitory effect on protests?, and
	\item  Why is it that some regimes fall in a matter of weeks, others fight to a stalemate, and still others survive relatively unscathed?
\end{enumerate}

We then established a simple one-compartment model that described the dynamics of a revolution sweeping through a population that is near-unanimous in its intrinsic dislike of the current dictatorial regime.  Despite the crudeness of our model, we are able to identify four main parameter regions that correspond to realistic situations in such countries: stable police state, meta-stable police state, unstable police state, and failed state.  These regions capture, at least qualitatively, a wide range of scenarios observed in the context of revolutionary movements in countries ruled by dictatorial regimes.  We examined two scenarios in detail, Tunisia and Egypt.  In the case of Tunisia we assumed that parameters were fixed for the lifetime of the revolution and explored the effect of different parameter regimes on the dynamics of the revolution.  In contrast, for Egypt we had to assume that parameters evolved with the revolution due to the importance of singular external events that were outside the scope of the basic model.  In both cases, our results qualitatively matched data describing frequency of Twitter and Blog posts by subject, as presented in \citet{HowardEtAl11}.  Finally, we briefly discussed how our model may further be able to describe aspects of the situation in 2009 Iran, and present-day China and Somalia.  We summarized these findings in Figure \ref{fig:Summary}.

\paragraph{}  From the Tunisian case study we concluded that the emergence of a small number of social media users and a low Internet penetration rate may have a large effect on the likelihood of a revolution.  This effect occurs because even a small increase in visibility or enthusiasm of revolutionaries can significantly decrease the basin of attraction of total state control.  This effect is further magnified by the nonlinear relationship between the magnitude and likelihood of shocks that a regime is subject to.  This phenomenon, together with the difficulty in establishing the parameters of models like ours, also explains why regimes manage to appear stable until the revolution is underway.  The Egyptian case study proposed a scenario which explains the limited impact of Internet shutdown on protests.  Essentially, by the time the Internet shutdown was fully implemented the revolution was sufficiently large that individuals' awareness was not impeded.  Finally, our review of Iran, China, and Somalia illustrates how our model is flexible enough to capture many different situations present in the world today.  This allows us to conclude that, despite the crudeness of our model, it is nevertheless valuable as a qualitative and conceptual tool.  It also indicates that our basic modelling assumptions appear to capture essential components of the dynamics of revolutions in dictatorial regimes.

\paragraph{}  The simplicity of our model also leaves much room for further work.  The adoption of a one-compartmental model has required the assumption of a homogeneous population.  One extension would be to expand the model to include additional compartments that take into account the heterogeneous nature of the population.  For example, Internet and social media use in the Arab world is highly concentrated among the youth \citep{HowardEtAl11}, thus youth communicate between each other differently than they communicate with their elders, or than elders communicate amongst themselves.  We could therefore extend our model by adding a compartment that explicitly takes into consideration the youth component of the revolutionary movement.  Similarly, we have assumed that the entire population desires regime change.  This may not be a good assumption in some cases.  For example, in Syria the government has deep rooted support from Alawite, Christian, and other minorities totalling approximately 25\% of the population \citep{Holliday11}.  So, our model could be amended to include a separate compartment to explicitly model populations that are likely to stay loyal to the regime until it is defeated.  Such a model may also be more suitable for describing past events in Libya and Bahrain.  Another approach to improving this model is to be more careful when choosing functional forms for the growth and decay terms in \eqref{eq:ODE}.  One might do this, for example, by studying the relationship between the micro-level threshold behaviour of individuals \citep{CentolaEtAl05, Kuran91} and the macro-level behaviour of the growth term.  Finally, since externalities, such as support or advice for protesters from foreign sources, are thought to be important \citep{HowardEtAl11, ZhuoEtAl11} it may be advisable to model these external influences explicitly.

\paragraph{}  Computational power has improved drastically over the past few decades to the point of enabling highly detailed and sophisticated models.  However, this is not to say that there is no longer a role for simple models.  Simple models have the advantages of relying on fewer assumptions about individual and communal behaviour.  They also admit a complete and rigorous mathematical analysis.  For simple models, it is frequently possible to characterize the entire parameter space and fully understand the behaviour of the model under all parameter regimes.  Despite their simplicity, these models nevertheless enable a conceptual and qualitative understanding of the phenomenon in question as long as they capture essential processes in the dynamics of the events under study.  Simple models, therefore, remain valuable for establishing a basic framework that can aid in understanding complex phenomena.

\newpage
\bibliographystyle{plainnat}
\bibliography{biblio}

\end{document}